\documentclass{article}

\usepackage{amsmath}
\usepackage{amssymb}
\usepackage{graphicx}

\evensidemargin \oddsidemargin
\numberwithin{table}{section}
\numberwithin{figure}{section}

\newtheorem{theorem}{Theorem}[section]
\newtheorem{lemma}[theorem]{Lemma}

\def\F{\mathcal F}
\def\G{\mathcal G}
\def\H{\mathcal H}

\def\O{\mathcal O}
\def\R{\mathcal R}

\begin{document}

\begin{center}
{\large\bf A Discontinuous Galerkin - Front Tracking Scheme \\ and its Optimal$^2$ Error Estimation  }

\vskip.20in

Tong Sun   and  Adamou Fode   %%$^{1}$
%\\[2mm]

%%$^{1}$
Department of Mathematics and Statistics\\
Bowling Green State University \\
Bowling Green, OH 43403

\end{center}

\date{today}

%{\footnotesize

%\begin{article}
%\keywords{Long-time error estimation, stability-smoothing indicator, parabolic PDE, %numerical solution, subdomain, partially implicit scheme.}

%{\bf AMS subject class.} Primary: 65M15  Secondary: 65L20

%\authorrunninghead{SUN}
%\titlerunninghead{Consistency + Numerical Smoothing $\Rightarrow$ Convergence}

\setcounter{page}{1}   %% This command is optional.
                       %% May set page number only for first page in
                       %% issue, if desired.

%% &lt;&lt;== End of commands to be entered at Wiley

%%  Authors, start here ==&gt;&gt;

%\author{Tong Sun}
%\affil{Department of Mathematical and Statistics\\
%Bowling Green State University\\
%Bowling Green, Ohio, USA}

\abstract{ In \cite{Rumsey} and \cite{rumseyDissertation}, an error estimate of optimal convergence rates and optimal error propagation (optimal$^2$) was given for the Runge-Kutta discontinuous Galerkin (RKDG) method solving the scalar nonlinear conservation laws in the case of smooth solutions. This manuscript generalizes the problem to the case of a piecewise smooth solution containing one fully developed shock. A front tracking technique is incorporated in the RKDG scheme to produce a numerical solution with a truly high order error. The numerical smoothness approach of \cite{Rumsey} is generalized to this particular case of a discontinuous solution.   }

\section{Introduction}

One of the most popular numerical methods for solving the scalar nonlinear conservation laws is the combination of the discontinuous Galerkin finite elements and the third order total variation diminishing Runge-Kutta scheme (RKDG). The first high order error analysis for the method was given by Zhang and Shu in \cite{qz2009}. The $L^2$ error estimate in \cite{qz2009} is of optimal convergence rates in space and time. Reading the proofs of \cite{qz2009}, one can see that the focus is on how the total variation diminishing Runge-Kutta scheme of order 3 (TVDRK-3) propagates error. In \cite{Rumsey}, since the smoothness of a numerical solution is employed for local error analysis, the $L^1$ contraction property of the scalar nonlinear conservation laws can be directly used for error propagation. Consequently, the $L^1$ error estimate of \cite{Rumsey} is not only of optimal convergence rates in space and time, but also of optimal error propagation, namely simple error accumulation. For brevity, we refer this kind of error estimates to be optimal$^2$ (read as: optimal square).

This manuscript shall generalize the result of \cite{Rumsey} to a special case of a discontinuous solution. Namely, the solution is supposed to consist of two smooth pieces separated by a shock. The Rankine-Hugoniot shock condition is employed as a differential equation to track the shock location in time. The discontinuous Galerkin method is used for the rest of the domain, with the shock as a moving interface and cell boundary point. The Rankine-Hugoniot equation and the DG scheme are combined to form one ODE system, the semi-discrete scheme. The fully discrete scheme is to solve the semi-discrete system by TVDRK-3. The semi-discrete scheme and the fully discrete scheme are both of optimal convergence order. An error estimate similar to the one in \cite{Rumsey} is proven, where the $L^1$ norm error grows linearly (simple accumulation). To the best knowledge of ours, there has not been such an error analysis for a shock solution in the literature.

Between the results of \cite{Rumsey} on smooth solutions and the results here in this paper on a fully developed shock, obviously one needs a method and its error analysis for solving the nonlinear conservation laws during shock formation. In this sense, if \cite{Rumsey} is considered to be the first step (smooth solution), this paper should be considered as the second step (simple fully developed shock). It is necessary to fill in the gap (transition from smooth solution to shock) in the future work. Beyond the above, a lot of work will be needed for the cases of contact discontinuities of zero or higher order, interaction between shocks, two and higher dimensional conservation laws, systems, and so on. Moreover, as a popular opinion, it is more interesting to estimate error for shock-capturing methods than for shock-tracking methods. Nevertheless, the result of this paper has its own significance. In case it is feasible and desirable to compute a shock solution very accurately, this paper tells us that it is possible to compute a shock with a truly high order approximation, and the error can be sharply estimated.  

As a matter of fact, the proofs for the case of smooth solutions in \cite{Rumsey}  apply to the smooth pieces of the solutions in this paper almost {\it verbatim}. Therefore only the proofs related to the shock and its neighboring cells shall be given in this paper. Results of numerical experiments will be given to show the boundedness of the smoothness indicators and that, as predicted by the error analysis, the computed shock location converges at the rate of $h^4$. Many details  of the lengthy error analysis for the smooth pieces of the solutions can be found in Rumsey's published Ph.D. dissertation \cite{rumseyDissertation}. For the new proofs, we will also show enough details to convince the reader on the error estimates and the boundedness and the computability of the coefficients, but leave the tedious details to a technical report (to be made public later).

The rest of the paper is organized as follows. In section 2, we present the conservation law and an example of a solution with a single fully developed shock. In section 3, we describe the fully discrete DG-FT-RK scheme. In section 4, we show the error analysis framework, which leads to the definition of the numerical smoothness indicators. In section 5, we state and prove the optimal$^2$ error estimate. In section 6, we show the numerical experiments, including the smoothness indicators and their boundedness, the anti-smoothing behavior of the scheme under inappropriate choices of time step size, and the optimal convergence rate of the shock location.

\section{The conservation law and the shock}

Consider the one-dimensional nonlinear conservation
law
\begin{equation}
u_t + f(u)_x = 0 \label{law}
\end{equation}
in a bounded interval $\Omega = [a,b]$.  In order to focus on the new ideas and the
new tools of the proof, we stay with the simple case of west wind ( $f'(u) > 0$ ) . 
Let the initial condition be
\begin{equation} \label{initialU}
u(0,x) = u_I(x)
\end{equation}
and the upwind boundary condition be
\begin{equation} \label{boundU}
u(t,a) = u_a(t).
\end{equation}
Assume that $u_I(x)$ is in $W^{p+1}_{\infty}(a,c)$ and $W^{p+1}_{\infty}(c,b)$,  except for that it contains a discontinuity at $x_s(0)=c \in (a,b)$, which is the initial position of a shock $x_s(t)$ for $t\ge 0$. 
Assume that the flux function $f(u)$ is sufficiently smooth and the
initial and boundary conditions are smooth and consistent to
guarantee that the entropy solution $u(t,x)$ remains in  $W^{p+1}_{\infty}(a,x_s(t))$ and $W^{p+1}_{\infty}(x_s(t),b)$ 
for all $t \in [0,T]$. 

The shock location is determined by the Rankine-Hugoniot shock condition, which is presented in the form of the following differential equation
\begin{equation} \label{Rankine}
\frac{d x_s}{d t} = \frac{f(u(x_s^-))-f(u(x_s^+))}{u(x_s^-)-u(x_s^+)} .
\end{equation}

\begin{figure}
  \includegraphics[width=6.2in,height=3.6in]{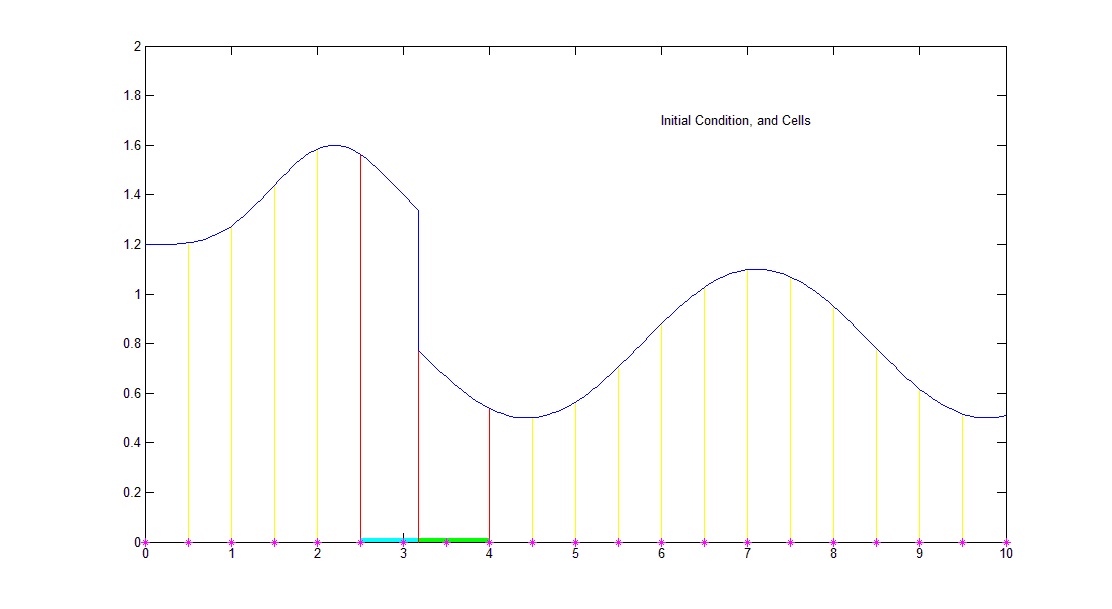}
  \caption{Initial condition and initial partition of a fully developed shock}\label{initialShock}
\end{figure}
An example of an initial condition with a shock is plotted in Figure \ref{initialShock}, where $(a,b) = (0,10)$, $x_s(0) = c = 3.18$, and
\begin{equation*}
u(0,x) = \begin{cases} 1.2 + 0.4 \, \sin^4 (x/1.4),              &    \text{if $ x \in (0,x_s(0))$}, \\
                                      0.8 - 0.3 \, \sin((x-3.1)/0.85),          &  \text{if $ x \in (x_s(0), 10)$}. \end{cases}
\end{equation*}
The boundary condition $u_a(t)=1.2$ is consistent with the initial condition within $W^{4}_{\infty}(0,x_s(0))$. This set of initial and boundary conditions will be used for the numerical experiments later in this paper, where we solve the Burgers' equation, $f(u)=u^2/2$ . 

It is easy to see that the entropy solution of the example will show most of the typical phenomena around a shock, including the solution's further sharpening, the growth of the shock height, the extrema of the smooth pieces to be kept constant before swallowed by the shock, the subsequent total variation diminishing, etc..    We have chosen this example to display all these phenomena in the numerical solution. The numerical solution contains neither smearing nor oscillation, which is necessary for the error to be of high order.

\section{The DG-FT-RK scheme}
The scalar nonlinear conservation law and the RKDG schemes are 
well-known \cite{ShuD}. The mass-conservative front tracking techniques can be found in \cite{Glimm} and the references there-in. We use a DG scheme with piecewise polynomials of degree $p$. Since the strength of the DG scheme is to obtain high order convergence, $p$ is supposed to be an integer bigger than or equal to $2$. All the numerical results in this paper are for $p=3$.  The Rankine-Hugoniot equation will be incorporated into the semi-discrete DG scheme to form a system of ordinary differential equations (ODE), where $x_s(t)$ is one of the state variables. Then the entire system is solved by the TVDRK-3 scheme to achieve a truly high order approximation to the exact solution of the conservation law, including the shock location. This section discribes the discontinuous Galerkin - front tracking - Runge-Kutta  scheme (DG-FT-RK). 

First, partition $\Omega$ with $a=x_{-1/2} < x_{1/2} < \cdots < x_{m-1/2} =
b$. Let $h = x_{j+1/2} - x_{j-1/2}$ be the same for all cells $\Omega_j = ( x_{j-1/2} , x_{j+1/2} )$. In the time direction, set the partition by $t_n=n \tau$, where $\tau$ is the uniform time step size. The uniform cell size and time step size are just for conveniences, not essential.

At time $t_n$, denote the computed shock location by $x_{sc,n}$. If $x_{sc,n} \in [x_{i-1/2},x_{i+1/2})$, then the three regular cells $\Omega_{i-1}$, $\Omega_{i}$ and $\Omega_{i+1}$ are converted to two special cells $\Omega_L = ( x_{i-3/2},x_{sc,n} )$ and $\Omega_R =  ( x_{sc,n},x_{i+3/2} )$. Figure \ref{initialShock} shows  the initial condition of the solution, and all regular cells and the two special cells at the initial time $t_0$ for the example, where $m=20$.

We denote the computed numerical solution by the pair $(x_{sc,n},u_n^c(x))$ at time $t_n$.  At the initial time $t_0=0$, $x_{sc,0}=x_s(0)$, $u_0^c(x)$ is taken to be the $L^2$-projection of $u_I(x)$ into piecewise discontinuous polynomials in the regular and special cells (given above), see Figure \ref{initialShock}. The description of the numerical scheme is inductive. Namely, we just show how to compute $(x_{sc,n+1},u_{n+1}^c(x))$ from $(x_{sc,n},u_{n}^c(x))$.

\subsection{The semi-discrete scheme}

Instead of having one semi-discrete solution satisfying the initial condition $u^h(0,x)=u_0^c(x)$, we employ a new semi-discrete solution $(x_{sh}(t),u^h(t,x))$ for each Runge-Kutta time step $[t_n,t_{n+1}]$, which satisfies the initial condition  
$$
x_{sh}(t_n) = x_{sc,n},  \qquad  \qquad   u^h(t_n,x)= u_n^c(x).
$$
Since each semi-discrete solution only lives in one time step, we use the same notation $u^h$ for all of them without worrying about the ambiguity. In fact, we always focus on one time step at a time. A part of a semi-discrete solution is its shock location $x_{sh}(t)$. Corresponding to the semi-discrete solution, the two special cells near the shock also evolve with the shock location. Namely, 
$$ \Omega^h_{L} = [x_{i-3/2},x_{sh}(t)], \qquad \qquad \Omega^h_{R} =  [x_{sh}(t),x_{i+3/2}].$$

To solve the problem with the
discontinuous Galerkin method, in each regular cell away from the shock and the two special cells surrounding the shock, we approximate the solution by using a polynomial in $\Pi_p$,  where $\Pi_p$ is the set of all the polynomials of
degree less than or equal to $p$. In each regular cell $\Omega_j$, $j = 0,\cdots,i-2$ and $i+2,\cdots,m-1$, a
semi-discrete solution $u^h$ satisfies, for any $v \in\Pi_p $,
\begin{equation}
(u^h_{\, t}, v)_{\Omega_j} = (f(u^h), v_x)_{\Omega_j}
                        + f(u^h(x_{j-1/2}^-)) v(x_{j-1/2}^+)
                        -  f(u^h(x_{j+1/2}^-)) v(x_{j+1/2}^-). \label{semi}
\end{equation}
Here the Godunov flux is employed under the west wind assumption
(for simplicity). At the upwind boundary $x_{-1/2}=a$, we set
\begin{equation} \label{boundUh}
u^h(t,x^-_{-1/2}) = u(t,a) = u_a(t).
\end{equation}
For the two special cells $\Omega^h_{L}$ and $\Omega^h_{R}$, the DG scheme is applied with a moving cell boundary at the semi-discrete shock location $x_{sh}(t)$. The formulations are
\begin{equation}
(u^h_{\, t}, v)_{\Omega^h_{L}} = (f(u^h), v_x)_{\Omega^h_{L}}
                        + f(u^h(x_{i-3/2}^-)) v(x_{i-3/2}^+)
                        - f(u^h(x_{sh}^-)) v(x_{sh}^-). \label{semiL}
\end{equation}
and
\begin{equation}
(u^h_{\, t}, v)_{\Omega^h_{R}} = (f(u^h), v_x)_{\Omega^h_{R}}
                        + f(u^h(x_{sh}^+)) v(x_{sh}^+)
                        - f(u^h(x_{i+3/2}^-)) v(x_{i+3/2}^-). \label{semiR}
\end{equation}
The semi-discrete shock location $x_{sh}(t)$ satisfies $x_{sh}(t_n) = x_{sc,n}$ and
\begin{equation} \label{semiRankine}
\frac{d x_{sh}}{d t} = \frac{f(u^h(x_{sh}^-))-f(u^h(x_{sh}^+))}{u^h(x_{sh}^-)-u^h(x_{sh}^+)} .
\end{equation}
The entire semi-discrete scheme consists of (\ref{semi}), (\ref{semiL}), (\ref{semiR}) and (\ref{semiRankine}). Since the time step size $\tau$ has to be sufficiently small and $x_{sh}(t)$ in the time step $[t_n,t_{n+1}]$ will never pass $x_{i+1/2} + h/2$ in practice, we have no need to consider redoing the special cells for a semi-discrete solution.

\subsection{The TVDRK-3 scheme}

The fully discrete scheme is essentially to use the TVDRK-3 scheme in each time step on the semi-discrete system (\ref{semi}, \ref{semiL}, \ref{semiR}, \ref{semiRankine}) of this time step. Abstractly, the TVDRK-3 scheme \cite{ShuD} for an ODE system $ W_t = \Psi(W) $ is given by 
\begin{eqnarray} \label{stage1}
   W^{1}_n  &&= W_n + \tau_n \Psi(W_n), \\
 \label{stage2}
 W^{2}_n &&= \frac{3}{4} W_n +  \frac{1}{4} W^{1}_n  + \frac{\tau_n}{4} \Psi(W^{1}_n), \\
 \label{stage3}
 W_{n+1}  &&= \frac{1}{3} W_n +  \frac{2}{3} W^{2}_n  + \frac{2 \tau_n}{3} \Psi(W^{2}_n).
\end{eqnarray}
Here we omit the detailed substitution of the semi-discrete system into this TVDRK-3 scheme because it is a standard procedure. Anyway, the TVDRK-3 scheme will compute $x_{sc,n+1}$ and $\hat{u}_{n+1}^c(x)$ for us from  $(x_{sc}(t_{n}), u_{n}^c(x))$ . 
$\hat{u}_{n+1}^c(x)$ is almost $u_{n+1}^c(x)$, but we still need the following.

\subsection{Replacing $\Omega_L$ and $\Omega_R$, if needed}

The only remaining task of the fully discrete scheme is to redo the two special cells when it is needed. For most of the time steps, if $x_{sc,n} \in [x_{i-1/2},x_{i+1/2})$, we shall find $x_{sc,n+1}$ staying in the interval. In this case, the $\hat{u}_{n+1}^c(x)$ computed by the TVDRK-3 is accepted as $u_{n+1}^c(x)$, and the work of the time step $[t_n,t_{n+1}]$ is completed. 

If  $x_{sc,n+1} \ge x_{i+1/2}$, it must be in $ [x_{i+1/2},x_{i+3/2}) = \{x_{i+1/2}\} \cup \Omega_{i+1} $. Now it is necessary to redo $\Omega_L$ and $\Omega_R$, and convert $\hat{u}_{n+1}^c(x)$ into a piecewise polynomial in the new cells. The result of this conversion will be accepted as $u_{n+1}^c(x)$ . The conversion of the cells and the solution is very simple. 

\begin{enumerate}

\item For the old $\Omega_L$, due to the moving of the shock location, now it has become $(x_{i-3/2},x_{sc,n+1})$. We split it into the regular cell $\Omega_{i-1} = (x_{i-3/2},x_{i-1/2})$ and the new $\Omega_L = (x_{i-1/2},x_{sc,n+1})$. Since $\hat{u}_{n+1}^c(x)$ restricted in the old $\Omega_L$ is a polynomial, we can simply take it as $u_{n+1}^c(x)$ in both  $\Omega_{i-1}$ and the new $\Omega_L$, hence no error is made here.  However, the polynomial needs to be represented under the bases of the new cells in the implementation.

\item For the old $\Omega_R$, now $(x_{sc,n+1},x_{i+3/2})$, we combine it with $\Omega_{i+2} = (x_{i+3/2},x_{i+5/2})$ to form the new $\Omega_R=(x_{sc,n+1},x_{i+5/2})$. In the new $\Omega_R$, $u_{n+1}^c(x)$ is taken to be the $L^2$-projection of $\hat{u}_{n+1}^c(x)$. The error of the projection will be estimated.

\item All the other regular cells remain unchanged. In each of them,  $\hat{u}_{n+1}^c(x)$ is accepted as  $u_{n+1}^c(x)$.

\end{enumerate}
Now the description of the fully discrete scheme is completed. It is obviously mass-conservative.

\section{Error analysis framework and smoothness indicators}

Since the error analysis for the smooth solutions in \cite{Rumsey} applies to the smooth pieces of the solution here, in the following sections we refer many details to \cite{Rumsey}. For the convenience of the reader, we stay with the same notations as much as possible. 

The goal of the error analysis is to estimate the error between the true solution $(x_s(t_n),u(t_n,x))$ and the numerical solution $(x_{sc,n}, u_n^c(x))$ at any time $t_n$. While $x_{sc,n}$ is the shock location of $u_n^c(x)$, the other discontinuities of $u_n^c(x)$ are just technical ones due to the DG approximation.

For each time step $[t_n,t_{n+1}]$, $n=0,1,\cdots, N-1$, an entropy solution $\tilde{u}(t,x)$ is used as an auxiliary function. We only need to specify the initial and boundary condition, so that $\tilde{u}(t,x)$ is uniquely determined as an entropy solution of (\ref{law}). Namely, the initial condition is $\tilde{u}(t_n,x) = u_n^c(x)$, and the boundary condition is $\tilde{u}(t,a)=u_a(t)$. Obviously, because $u_n^c(x)$ is a DG function, there is a technical discontinuity at each cell boundary point $x_{j-1/2}$, for $j=1,\cdots,i-1$ and $i+2,\cdots,m-1$. Since $\tilde{u}(t,x)$ evolves out of $u_n^c(x)$, it contains many discontinuities. Nevertheless, we define 
$$ \|\tilde{u}(t_{n+1},x) - u_{n+1}^c(x) \|_{L^1(\Omega)} $$
as the local error. Consequently, we have the optimal error propagation estimate 
$$
\|u(t_{n+1},x) - \tilde{u}(t_{n+1},x) \|_{L^1(\Omega)} \le \|u(t_{n},x) - \tilde{u}(t_n,x) \|_{L^1(\Omega)}  = \|u(t_{n},x) - u_n^c(x) \|_{L^1(\Omega)}
$$
due to the well-known $L^1$-contraction between entropy solutions, see Chapter 16 of \cite{smoller}. Now, the macro framework of the error analysis is
\begin{equation} \label{propagation}
\|u(t_{n+1},x) -u_{n+1}^c(x) \|_{L^1(\Omega)} \le \|u(t_{n},x) - u_n^c(x) \|_{L^1(\Omega)} + \|\tilde{u}(t_{n+1},x) - u_{n+1}^c(x) \|_{L^1(\Omega)}.
\end{equation}

While we are pleased with the optimal and trivial error propagation, we are paying the price of a non-standard and highly non-trivial analysis on the local error  $ \|\tilde{u}(t_{n+1},x) - u_{n+1}^c(x) \|_{L^1(\Omega)} $. We will further split the local error into several parts, and estimate each part carefully to achieve optimal convergence rates in space and time. The fundamental reason that the local error can be estimated with optimal convergence rates  is the {\it numerical smoothness} of the computed solution $u_n^c(x)$ and a few auxiliary functions (all closely related to $u_n^c(x)$ ). Therefore we have to start from a set of numerical smoothness indicators, which only depends on $u_n^c(x)$. %\section{The smoothness indicators}
We define the following spatial and temporal smoothness indicators for each time step $[t_n,t_{n+1}]$.
\begin{itemize}
\item Spatial smoothness indicator: $S_n^{p} = S_{p}(u^c_n),$
\item Temporal smoothness indicator: $T_n^{k} = T_{k}(u^c_n).$
\end{itemize}
Here $S$ stands for space, $T$ stands for time, $p$ is the degree of
the polynomials in each cell, and $k$ is the order of the
Runge-Kutta scheme.

\subsection{Definition of $S^p_n$}

The spatial smoothness indicator $S^p_n$ contains the
spatial derivatives of $u^c_n$ within each cell, as well as the modified measurements of the jumps
of the derivatives across the regular cell boundaries. However, there is no jump information needed at the shock, which is the moving boundary point between the two special cells $\Omega_L$ and $\Omega_R$. In fact, the approximation accuracy of the shock only depends on the smoothness of the solution in the two special cells. 

At each regular cell boundary point $x_{j-1/2}$, $j=0, 1,\cdots,i-1$ and $i+2,\cdots,m-1$, let us consider the right limits of the derivatives of $u_n^c(x)$ of order $l=0,1,\cdots,p$,
$$ 
M_{n,j}^l = \frac{\partial^l}{\partial x^l} u^c_n(x_{j-1/2}^+).
$$
The values of these $M_{n,j}^l$ deliver the information about the smoothness of $u_n^c(x)$ in the interior of each regular cell. In the case $j=i-1$, the values of $M_{n,i-1}^l$ actually tells us the smoothness of $u_n^c(x)$ in the interior of the special cell $\Omega_L$. For the other special cell $\Omega_R$, the smoothness of $u_n^c(x)$ in $\Omega_R$ can be delivered by
$$ 
M_{n,s}^l = \frac{\partial^l}{\partial x^l} u^c_n(x_{sc,n}^+).
$$

In order to extract the information about the numerical smoothness of $u_n^c(x)$ at each regular cell boundary point $x_{j-1/2}$, we also compute the left limits of the derivatives of $u_n^c(x)$,
$$ 
L_{n,j}^l = \frac{\partial^l}{\partial x^l} u^c_n(x_{j-1/2}^-).
$$
For the case $j=0$, see \cite{Rumsey} for the details. Now we can compute the jumps
$$
J_{n,j}^l = M_{n,j}^l - L_{n,j}^l 
$$
and, for $\alpha = 1/p$,
$$
D_{n,j}^l = \frac{J_{n,j}^l}{h^{p+2 -l(1+\alpha)}},
$$
for all $l=0,1,\cdots,p$ and $j=0,1,\cdots,i-1, i+2,\cdots,m-1$. Since the values of $L_{n,j}^l$ can be computed from the values of $M_{n,j-1}^l$ by Taylor expansions, we have included redundant information already. All the smoothness information is contained in the $M$'s, but the $D$'s deliver the jump smoothness information explicitly.

Formally, the smoothness indicator $S_n^p$ is 
$$
S_n^p = \left( M_{n,s}^l, M_{n,j}^l, D_{n,j}^l \, | \, l=0,1,\cdots,p, j=0,1,\cdots,i-1, i+2,\cdots,m-1 \right).
$$

It is obvious that the values of $M^l_{n,j}$ and $M^l_{n,s}$ should
be of $\O(1)$, which means that $u_n^c(x)$ is smooth in the interior of each regular and special cell. A necessary condition has been given in \cite{necessity} on the jumps
$J_{n,j}^l$. If $u_n^c(x)$ is approximating $u(t_n,x)$ in its smooth pieces in the optimal rate $h^{p+1}$, the jumps must satisfy
$$
|J_{n,j}^l | \le C h^{p+1-l}.
$$
In the case of the DG method, \cite{Rumsey} and here, both error analysis and
numerical experiments suggest that the jumps are even smaller. In fact, $D_{n,j}^l =
J_{n,j}^l/h^{p+2-l(1+\alpha)}$ should be at most of $\O(1)$ for $j=0,\cdots,i-1,i+2,\cdots,m-1$. In the error estimates, $D_{n,j}^l$ play the role of high order derivatives as in the traditional error analysis. Therefore, it is natural to have $D_{n,j}^l$ serving as a part of the smoothness indicator.

\subsection{Definition of $T^k_n$}

The temporal smoothness indicator $T_n^{k}$ consists of the temporal
derivatives of $x_{sh}(t)$ and $u^{h}(t,x)$ at $t=t_n$. Namely,
$$
 T_n^{k} = \left( \left[ x_{sc,n},u^c_n \right], \left[\frac{d}{dt}x_{sh}(t_n),u^{h}_{t}(t_n)\right], \left[\frac{d^2}{dt^2}x_{sh}(t_n),u^{h}_{tt}(t_n)\right],
\cdots, \left[ \frac{d^{k+1}}{dt^{k+1}} x_{sh}(t_n), \frac{\partial^{k+1} }{\partial t^{k+1}}u^{h}(t_n) \right] \right).
$$ 
The first derivatives $\frac{d}{dt}x_{sh}(t_n)$ and  $u^{h}_{t}(t_n)$ is computed by using (\ref{semi},\ref{semiL},\ref{semiR},\ref{semiRankine}).

The formula for computing the second derivative can be obtained by
taking derivatives with respect to $t$ on both sides of the
semi-discrete equations (\ref{semi},\ref{semiL},\ref{semiR},\ref{semiRankine}). For the equation (\ref{semi}) about the regular cells, it is easy to derive
\begin{eqnarray}
(u^h_{tt}, v)_{\Omega_j} &=& (f'(u^h) u^h_{t}, v_x)_{\Omega_j} + f'(u^h(x_{j-1/2}^-)) u^h_{t}(x_{j-1/2}^-) v(x_{j-1/2}^+) \nonumber \\
        && - f'(u^h(x_{j+1/2}^-)) u^h_{t}(x_{j+1/2}^-) v(x_{j+1/2}^-). \label{secondD}
\end{eqnarray}
For the special cells $\Omega_{L}^h$ and $\Omega_{R}^h$, we take derivative with respect to $t$ on both sides of (\ref{semiL}) and (\ref{semiR}) to compute $u^h_{tt}(t_n,x)$. After the cancellations and regrouping of some terms, we get
\begin{eqnarray}
(u^h_{tt}, v)_{\Omega_L^h} &=& (f'(u^h) u^h_{t}, v_x)_{\Omega_L^h}  
+ f'(u^h(x_{i-3/2}^-)) u^h_{t}(x_{i-3/2/2}^-) v(x_{i-3/2}^+)  \nonumber \\
&-& f'(u^h(x_{sh}^-)) u^h_{t}(x_{sh}^-) v(x_{sh}^-) \nonumber \\ 
&-& \left[ u^h_t(x_{sh}^-) + f'(u^h(x_{sh}^-)) u^h_x(x_{sh}^-) \right] \, \frac{d x_{sh}}{dt} \,\,\,  v(x_{sh}^-) , \label{secondL}
\end{eqnarray}
and
\begin{eqnarray}
(u^h_{tt}, v)_{\Omega_R^h} &=& (f'(u^h) u^h_{t}, v_x)_{\Omega_R^h}  
+ f'(u^h(x_{sh}^+)) u^h_{t}(x_{sh}^+)  v(x_{sh}^+)  \nonumber \\
&-& f'(u^h(x_{i+3/2}^-)) u^h_{t}(x_{i+3/2/2}^-) v(x_{i+3/2}^-) \nonumber \\
&+& \left[ u^h_t(x_{sh}^+) + f'(u^h(x_{sh}^+)) u^h_x(x_{sh}^+) \right] \, \frac{d x_{sh}}{dt} \,\,\, v(x_{sh}^+)  . \label{secondR}
\end{eqnarray}
Differentiating (\ref{semiRankine}) with respect to $t$ by the product rule,
\begin{eqnarray}
\frac{d^2 x_{sh}}{dt^2} &=& \frac{f'(u^h(x_{sh}^-))\left[u^h_t(x_{sh}^-)+u^h_x(x_{sh}^-)\frac{d\, x_{sh}}{dt} \right]
-f'(u^h(x_{sh}^+))\left[u^h_t(x_{sh}^+)+u^h_x(x_{sh}^+)\frac{d\, x_{sh}}{dt} \right] }{u^h(x_{sh}^-) - u^h(x_{sh}^+) } \nonumber \\
&-&   [f(u^h(x_{sh}^-))-f(u^h(x_{sh}^+))] \frac{\left[u^h_t(x_{sh}^-)+u^h_x(x_{sh}^-)\frac{d\, x_{sh}}{dt} \right]-\left[u^h_t(x_{sh}^+)+u^h_x(x_{sh}^+)\frac{d\, x_{sh}}{dt} \right]}{ [u^h(x_{sh}^-) - u^h(x_{sh}^+) ]^2 }  \nonumber \\
&=& \frac{ \left[ f'(u^h(x_{sh}^-)) - \frac{d\, x_{sh}}{dt}\right] \left[u^h_t(x_{sh}^-)+u^h_x(x_{sh}^-)\frac{d\, x_{sh}}{dt} \right]}{u^h(x_{sh}^-) - u^h(x_{sh}^+) } \nonumber \\
&+&  \frac{\left[\frac{d\, x_{sh}}{dt} - f'(u^h(x_{sh}^+)) \right] \left[u^h_t(x_{sh}^+)+u^h_x(x_{sh}^+)\frac{d\, x_{sh}}{dt} \right] }{u^h(x_{sh}^-) - u^h(x_{sh}^+) }  .  \label{secondS}
\end{eqnarray}
The term $\left[\frac{d^2}{dt^2}x_{sh}(t_n),u^{h}_{tt}(t_n)\right]$ of $T_n^k$ can be computed by substituting $t=t_n$ in (\ref{secondD},\ref{secondL},\ref{secondR},\ref{secondS}). Although the formulas of the computation seem to be tedious, one can easily see the boundedness of this component of $T_n^k$ from the boundedness of the spatial smoothness indicator $S_n^p$. 

In the last term of the right hand side of (\ref{secondL}) and (\ref{secondR}), the factor $ u^h_t(x_{sh}^\pm) + f'(u^h(x_{sh}^\pm)) u^h_x(x_{sh}^\pm)$ should be nearly $0$, because their true solution counterpart $ u_t(x_{s}^\pm) + f'(u(x_{s}^\pm)) u_x(x_{s}^\pm) $ equals to $0$. Therefore, in practice, $u^h_{tt}(t_n,x)$ in $\Omega^h_L$ and $\Omega^h_R$ can be approximated, with the last term of (\ref{secondL}) and (\ref{secondR}) neglected.

The later components of $T_n^k$ can be derived and computed in the same way. We omit the details because of the increasing tediousness of the formulas. We just emphasize that it is easy to see that all the components of $T_n^k$ are computable, and their boundedness follows the boundedness of the spatial smoothness indicator. Similar to the case of $u^h_{tt}$, we can neglect some terms in the computation of  $u^h_{ttt}$ and $u^h_{tttt}$ in $\Omega^h_L$ and $\Omega^h_R$, leaving the cost of computing $u^h_{ttt}(t_n,x)$ and $u^h_{tttt}(t_n,x)$ in the special cells the same as in the regular cells. 

{\bf Remark}: Since the boundedness of $T_n^k$ can be derived from the boundedness of $S_n^p$, the entire temporal smoothness indicator $T_n^k$  can be considered as redundant information in certain sense. However, the direct computation of $T_n^k$ can provide a sharper {\it a posteriori}  error estimate on the RK scheme. Moreover, it is observed in our numerical experiments that the high order discontinuities of a solution can be detected by the temporal smoothness indicator $T^k_n$ more sensitively (to be included in \cite{fodeDissertation}).

\section{The error analysiss}

\begin{theorem} \label{main}
Let $u(t,x)$ be the entropy solution of the nonlinear conservation
law (\ref{law}) satisfying the initial condition (\ref{initialU})
and upwind boundary condition (\ref{boundU}). Let $x_s(t)$ be the shock location of $u(t,x)$. Let $u^c_n$ be the
numerical solution computed by the DG-FT-RK scheme described in Section 3. Assume that $u$ and
$u^c_n$ are bounded by a constant $U$ in $[0,T] \times \Omega$. Let
$\beta = \max_{|w| \le U} f'(w)$. Assume that the time step size
$\tau$ satisfies the standard CFL condition $\beta
\tau \le h$ and the strengthened CFL condition $\tau \le \gamma
h^{1+\alpha}$, for a positive constant
$\gamma$, and $\alpha = 1/p$. Let $\tilde{u}(t,x)$ be as defined in Section 4. 

If there is a positive real number $M$, such that, for all $t_n \le
T$, all the components of $S^p_n$ and $T^k_n$ are bounded by $M$,
then the spatial and temporal local error in $[t_n,t_{n+1}]$ satisfy
\begin{equation}
\|\tilde{u}(t_{n+1},x) - u^{h}(t_{n+1},x)\|_{L_1(\Omega)} \le \tau
h^{p+1}  \F (S_n^{p}), \label{spaceError}
\end{equation}
\begin{equation}
\|u^{h}(t_{n+1},x) - \hat{u}^c_{n+1}(x) \|_{L_1(\Omega)} \le \tau^{k+1}
\G(T_n^{k},S_n^{p}), \label{timeError}
\end{equation}
where $\F(S^p_n)$ and $\G(T^k_n,S_n^{p})$ are computable functions
of the indicators. For the transition error occurring when $\Omega_R$ is redone, we define 
$$
h_{cr}(n) = \begin{cases} 0, & \Omega_R \text{ is not redone in the time step } [t_n,t_{n+1}], \\
h, &   \Omega_R \text{ is redone  in the time step } [t_n,t_{n+1}].
\end{cases}
$$
With this function defined, for a computable function $\mathcal{H}(\cdot)$ and $\hat{S}^p_{n+1}$ defined in Lemma \ref{projectionErrorOmegaR},
\begin{equation} \label{errorOmegaR} 
\| \hat{u}^c_{n+1}(x) - u^c_{n+1}(x) \|_{L_1(\Omega)} \le h_{cr}(n) h^{p+1} \mathcal{H}(\hat{S}^p_{n+1}).
\end{equation}

As a consequence of the the local error estimates (\ref{spaceError}),
(\ref{timeError}) and (\ref{errorOmegaR}), and the error splitting and the $L_1$-contraction property (\ref{propagation}),
$$
\|u(t_{n+1}) - u^c_{n+1}\|_{L_1(\Omega)} \le \|u(t_{n}) -
u^c_{n}\|_{L_1(\Omega)} + \tau [ h^{p+1} \F(S_n^{p}) +  \tau^k
\G(T_n^{k},S_n^{p}) ] + h_{cr}(n) h^{p+1} \mathcal{H}(\hat{S}^p_{n+1}).
$$
Finally, at the end of the computation ($t_N = T$),
\begin{equation} \label{global}
\|u(T) - u^c_N\|_{L_1(\Omega)} \le \|u(0) - u^c_0\|_{L_1(\Omega)} +
\sum_{n=0}^{N-1} \tau [ h^{p+1} \F(S_n^{p}) +  \tau^k
\G(T_n^{k},S_n^{p}) ] + \sum_{n=0}^{N-1} h_{cr}(n) h^{p+1} \mathcal{H}(\hat{S}^p_{n+1}),
\end{equation}
and consequently, for the shock location, let $H(T) = \frac{3}{4} |u(T, x_s^-(T)) - u(T, x_s^+(T))|$, then
\begin{equation} \label{globalShock}
 H(T) \, |x_s(T)-x_{sc,N}| \le \|u(0) - u^c_0\|_{L_1(\Omega)} + \mathcal{O}(h^{p+1} + \tau^k).
\end{equation}

\end{theorem}

\noindent {\bf Proof.} It suffices to prove (\ref{spaceError}),
(\ref{timeError}) and (\ref{errorOmegaR}), which will be done in the following three subsections. It is easy to see that $\sum_{n=0}^{N-1} h_{cr}(n) \le |\Omega|$, then the rest of the results are obvious. \# \\

{\bf Remark.} If $k=p$, then $\tau = \mathcal{O}(h^{1+\alpha})$ with $\alpha = 1/p$ implies $\tau^k = \mathcal{O}(h^{p+1})$. In the literature, some proofs require $\alpha = 1/3$ \cite{qz2004}, some other proofs allow $\alpha = 0$ \cite{qz2009}. Here we require  $\alpha = 1/p$. In the numerical experiments, we use $k=p=3$. The numerical experiments show that, when  $\alpha = 1/3$, the shock location converges in the optimal rate $h^4$ in space,  and $\tau^3$ in time. When we chose  $\alpha < 1/3$, the optimal rate in space ($h^4$) was lost, although we still had convergence. 

\subsection{Estimating $\tilde{u}(t_{n+1}) - u^{h}(t_{n+1})$, proof
of (\ref{spaceError}) }

We begin with introducing an auxiliary piecewise PDE solution
$u^{e}(t,x)$ and its shock location $x_{se}(t)$. First define a local strong solution $u^{e}_j$ of a Cauchy problem of the
conservation law (\ref{law}) for each regular cell $\Omega_j$, $j=0,\cdots,i-2,i+2,\cdots,m-1$. The initial values of $u^e_j$ are given on the
line segment $\{t_n\} \times \overline{\Omega_{j-1} \cup \Omega_j}$ by
\begin{equation} \label{uejo}
u^{e}_j(t_n,x) = M_{n,j}^0 + M_{n,j}^1 (x-x_{j-1/2}) + \cdots +
\frac{M_{n,j}^p}{p!} (x-x_{j-1/2})^p .
\end{equation}
It is easy to see that, in $\Omega_j$, $u^{e}_j(t_n) = u^c_n$. Moreover, in
$\Omega_{j-1}$,
$$
u^{e}_j(t_n) =u^c_n + J_{n,j}^0 + J_{n,j}^1 (x-x_{j-1/2}) + \cdots +
\frac{J_{n,j}^p}{p!} (x-x_{j-1/2})^p,
$$
in all regular cells (except for $\Omega_0$, see \cite{Rumsey} for the detail of $\Omega_0$). As a strong solution of the Cauchy problem,  $u^{e}_j$ certainly exists in the
region
 $$\R_{n,j} = \{(t,\tilde{x})| t \in [t_n,t_{n+1}], x \in
\overline{\Omega_{j-1} \cup \Omega_j}, \tilde{x} = x + f'(u^{e}_j(t_n,x))(t-t_n), \tilde{x} \le x_{j+1/2} \}.
$$
 This is the trapezoidal region shown
below, which is covered by the characteristic lines originating from
$\overline{\Omega_{j-1} \cup \Omega_j}$. When $\tau$ satisfies the standard
CFL condition $\beta \tau \le h$, it is easy to see that $[t_n,
t_{n+1}] \times \Omega_j \subset \R_{n,j} \subset [t_n, t_{n+1}]
\times \overline{\Omega_{j-1} \cup \Omega_j}$.

%\begin{figure}
\setlength{\unitlength}{1.0mm}
\begin{picture}(70,13)(-20,15)
    \linethickness{0.5pt}
        \put(20, 0){\line(1,0){60}}
        \put(20,20){\line(1,0){60}}
        \put(30, 0){\line(0,1){20}}
        \put(50, 0){\line(0,1){20}}
        \put(70, 0){\line(0,1){20}}

  % \linethickness{1.4pt}
        \put(70, 0){\line(0,1){20}}
  %  \linethickness{1.3pt}
        \qbezier(30,0)(35,10)(40,20)

    \put(16,-1){\makebox(0,0)[b]{$t_n$}}
    \put(15,19){\makebox(0,0)[b]{$t_{n+1}$}}
    \put(30,-5){\makebox(0,0)[b]{$x_{j-\frac{3}{2}}$}}
    \put(50,-5){\makebox(0,0)[b]{$x_{j-\frac{1}{2}}$}}
    \put(70,-5){\makebox(0,0)[b]{$x_{j+\frac{1}{2}}$}}
    \put(56, 8){\makebox(0,0)[b]{$\R_{n,j}$}}
\end{picture}
%\end{figure}

\vspace{1.1in}

At the upwind boundary, let $u^{e}_{-1} = u$ for $x \in \Omega_{-1}
= (x_{-1/2}-h,x_{-1/2})$. Due to the smoothness of $u_a(t)$, one can
determine the value of $u$ in $\Omega_{-1}$ by tracing back (but not
computable).  

Let $u^e_L(t,x)$ be the strong solution of the Cauchy problem with the initial value  
$$
u^{e}_L(t_n,x) = M_{n,i-1}^0 + M_{n,i-1}^1 (x-x_{i-3/2}) + \cdots + \frac{M_{n,i-1}^p}{p!} (x-x_{i-3/2})^p, \qquad x \in [x_{i-5/2},x_{sc,n}].
$$
The region where the solution of the Cauchy problem exists is  
$$\R_{n,L} = \{(t,\tilde{x})| t \in [t_n,t_{n+1}], x \in \overline{\Omega_{i-2} \cup \Omega_L}, \tilde{x} = x + f'(u^{e}_L(t_n,x))(t-t_n) \}. 
$$

Let $u^e_R(t,x)$ be the strong solution of the Cauchy problem with initial value  
$$
u^{e}_R(t_n,x) = M_{n,s}^0 + M_{n,s}^1 (x-x_{sc,n}) + \cdots + \frac{M_{n,s}^p}{p!} (x-x_{sc,n})^p, \qquad x \in [x_{sc,n},x_{i+3/2}].
$$
The region where the solution of the Cauchy problem exists is  
$$\R_{n,R} = \{(t,\tilde{x})| t \in [t_n,t_{n+1}], x \in \overline{\Omega_R}, \tilde{x} = x + f'(u^{e}_R(t_n,x))(t-t_n), \tilde{x} \le x_{i+3/2} \}. 
$$
$\R_{n,L}$ and $\R_{n,R}$ are shown in the picture below. It is important to notice that $\R_{n,L}$ and $\R_{n,R}$  have an intersection $\Delta ABC$, because the wave speed on the upwind side of the shock is bigger than the wave speed on the downwind side. The curve $x_{se}(t)$ gives the location of the auxiliary shock, sitting in the interior of the intersection of $\R_{n,L}$ and $\R_{n,R}$. 

\setlength{\unitlength}{1.0mm}
\begin{picture}(70,15)(-5,15)
    \linethickness{0.5pt}
        \put(20, 0){\line(1,0){100}}
        \put(20,20){\line(1,0){100}}

        \put(30, 0){\line(0,1){20}}
        \put(50, 0){\line(0,1){20}}
        \put(70, 0){\line(0,1){20}}
        \put(90, 0){\line(0,1){20}}
        \put(110, 0){\line(0,1){20}}

   % \linethickness{1.4pt}
        \put(110, 0){\line(0,1){20}}
   % \linethickness{1.3pt}
        \qbezier(30,0)(35,10)(40,20)
        \qbezier(75,0)(81,10)(87,20)
        \qbezier(75,0)(77,10)(79,20)
        \qbezier(75,0)(79,12.5)(83.5,20) %shock

	\linethickness{0.1pt}
	\put(75,0){\line(2,-1){4}}
	\put(83.5,20){\line(0,1){3.5}}

    \put(16,-1){\makebox(0,0)[b]{$t_n$}}
    \put(15,19){\makebox(0,0)[b]{$t_{n+1}$}}
    \put(30,-5){\makebox(0,0)[b]{$x_{i-\frac{5}{2}}$}}
    \put(50,-5){\makebox(0,0)[b]{$x_{i-\frac{3}{2}}$}}
    \put(67,-5){\makebox(0,0)[b]{$x_{i-\frac{1}{2}}$}}
    \put(90,-5){\makebox(0,0)[b]{$x_{i+\frac{1}{2}}$}}
    \put(110,-5){\makebox(0,0)[b]{$x_{i+\frac{3}{2}}$}}
    
    \put(57, 8){\makebox(0,0)[b]{$\R_{n,L}$}}
    \put(95, 8){\makebox(0,0)[b]{$\R_{n,R}$}}

    \put(81,-4.7){\makebox(0,0)[b]{$x_{sc,n}$}}
    \put(85,24){\makebox(0,0)[b]{$x_{se}(t_{n+1})$}}

     \put(74.5,-3){\makebox(0,0)[b]{$A$}}
     \put(79,20.5){\makebox(0,0)[b]{$C$}}
     \put(87,20.5){\makebox(0,0)[b]{$B$}}

\end{picture}

\vspace{0.95in}

\noindent The auxiliary shock $x_{se}(t)$ is determined by $x_{se}(t_n) = x_{sc,n}$ and
\begin{equation} \label{auxiiRankine}
\frac{d x_{se}}{d t} = \frac{f(u^e_L(x_{se}))-f(u^e_R(x_{se}))}{u^e_L(x_{se})-u^e_R(x_{se})} .
\end{equation}
It is easy to show that $f'(u^e_L(x_{se})) > \frac{d x_{se}}{d t}  = \frac{f(u^e_L(x_{se}))-f(u^e_R(x_{se}))}{u^e_L(x_{se})-u^e_R(x_{se})} > f'(u^e_L(x_{se}))$ for a shock.

Now, we are ready to define the local piecewise PDE solution. Let $\Omega^e_L = (x_{i-3/2},x_{se}(t))$ and $\Omega^e_R = (x_{se}(t),x_{i+3/2})$.
\begin{equation}
u^{e}(t,x) = 
\begin{cases}
u^{e}_{j}(t,x),     & t \in [t_n,t_{n+1}], \, x \in \Omega_j, \,  j=0,\cdots,i-2,i+2,\cdots,m-1, \\
u^e_L(t,x),           & t \in [t_n,t_{n+1}], \, x \in \Omega^e_L, \\
u^e_R(t,x),           & t \in [t_n,t_{n+1}], \, x \in \Omega^e_R.
\end{cases}
\end{equation}
The shock location of $\tilde{u}(t,x)$  is also $x_{se}(t)$, because it is easy to prove that $\tilde{u}(t,x)$ equals to $u^e_L(t,x)$ in $\Omega^e_L$ and it equals to $u^e_R(t,x)$ in $\Omega^e_R$.  Although $u^e(t,x)$ is not a solution of the conservation law (\ref{law}) globally in $\Omega$, but it is a solution locally, especially near the shock. The shock location of $u^e(t,x)$ is also $x_{se}(t)$. Now we have the following estimate on $\tilde{u}(t_{n+1}) - u^e(t_{n+1})$ in $\Omega$.
\begin{lemma} \label{lemmaTE}
If $\beta \tau \le h$ and  $\tau \le \gamma  h^{1+\alpha}$, then there is a computable function $\tilde{C}(S^p_n)$ depending on the spatial smoothness indicator $S^p_n$, such that, for any sufficiently small $\tau$,
$$
\|\tilde{u}(t_n+\tau) - u^e(t_n+\tau)\|_{L_1(\Omega)} \le \tau
h^{p+1} \tilde{C}(S^p_n).
$$
\end{lemma}
\noindent {\bf Proof.} Lemma \ref{lemmaTE} is a generalization of Lemma 4.4 of \cite{Rumsey}. First notice that $\tilde{u}(t,x) = u^e(t,x)$ in $\Omega^e_R$, hence the error here is zero. For all the regular cells, the proof of Lemma 4.4 in \cite{Rumsey} can be applied here {\it verbatim}. Since the only error between $\tilde{u}(t,x)$ and $u^e(t,x)$ occurs in a small sub-interval at the very left end of each cell, including cell  $\Omega^e_L$, the proof of Lemma 4.4 of \cite{Rumsey} applies to $\Omega^e_L$ as well. So we refer the reader to \cite{Rumsey}. \# \\

%\vspace{0.05in}

It remains to estimate $u^e(t,x)-u^h(t,x)$, but we need an $L^2$- projection $u^p(t,x)$ of $u^e(t,x)$ cell by cell. To be specific, for each regular cell $\Omega_j$,
$$
(u^{p}_j,v)_{\Omega_j} = (u^{e}, v)_{\Omega_j}, \qquad \forall v \in \Pi_p.
$$
For the special cells $\Omega_L$ and $\Omega_R$, we actually need to define the projection in the following expanded and overlapping intervals. Consider the two characteristic lines in the above picture at $x_{sc,n}$, namely $x_{s+}(t) =  x_{sc,n} + f'(u^{e}_L(t_n,x_{sc,n}))(t-t_n)$ and $x_{s-}(t) =  x_{sc,n} + f'(u^{e}_R(t_n,x_{sc,n}))(t-t_n)$, which are the right boundary of $\R_{n,L}$ and the left boundary of $\R_{n,R}$ respectively. We define the following overlapping projections $u^p_L(t,x)$ in $\Omega_{LE}=(x_{i-3/2},x_{s+}(t))$ and $u^p_R(t,x)$ in $\Omega_{RE}=(x_{s-}(t),x_{i+3/2})$, where
$$
 (u^p_L(t),v)_{\Omega_{LE}} = (u^{e}_L(t), v)_{\Omega_{LE}}, \qquad \forall v \in \Pi_p,
$$
and
$$
 (u^p_R(t),v)_{\Omega_{RE}} = (u^{e}_R(t), v)_{\Omega_{RE}}, \qquad \forall v \in \Pi_p.
$$
Since $u^e_L(t_n,x)$, $u^e_R(t_n,x)$ and $u^{e}_j(t_n,x)$ are polynomials,  if the smoothness indicators are bounded, $u^e_L(t,x)$, $u^e_R(t,x)$ and $u^{e}_j(t,x)$ are smooth in their domains respectively. To reveal more details on their smoothness, and the error of the projections in each cell, we have the following two lemmas.

\begin{lemma} \label{lemmaSS}
Consider $\Omega_j$ for $j=0,\cdots,i-2,LE,RE,i+2,\cdots,m-1$. When $j=LE$, it can be understood as $j=i-1$ or $j=L$ according to the context. Similarly, when $j=RE$, it can be understood as $j=s$ or $j=R$ according to the context.

 There are constants $N^l_{n,j}$  ($l=0,1,\cdots, p+1,  j=0,\cdots,i-2,LE,RE,i+2,\cdots,m-1$), which depend on the flux
 function $f$ and can be computed from $M^0_{n,j}, M^1_{n,j}, \cdots, M^p_{n,j}$, such
 that, for any $t \in (t_n,t_{n+1}]$,
$$
\left\| \frac{\partial^l}{\partial x^l} u^{e}_{j}(t,x) \right\|_{L_{\infty}(\R_{n,j})} \le N^l_{n,j}, \qquad l=0,1,\cdots,p.
$$
Moreover,
\begin{equation} \label{pPlus1}
\left \| \frac{\partial^{p+1}}{\partial x^{p+1}} u^{e}_{j}(t,x) \right\|_{L_{\infty}(\Omega_j)} \le (t-t_n) N^{p+1}_{n,j} .
\end{equation}
\end{lemma}
\noindent {\bf Proof.}  The proof  of Lemma 4.2 in \cite{Rumsey} applies to this Lemma {\it verbatim}.    \#  \\

The $t-t_n$ on the right hand side of (\ref{pPlus1}) is crucial. See the remarks and proofs of \cite{Rumsey}. 
The result of Lemma \ref{lemmaSS} makes the following projection error estimation obvious. 

\begin{lemma} \label{lemmaEP}
Consider $\Omega_j$, for $j=0,\cdots,i-2,LE,RE,i+2,\cdots,m-1$. For sufficiently small $\tau = t-t_n$,
\begin{equation} \label{EPL1}
\|u^e_j - u^p_j\|_{L_1(\Omega_j)} \le C_1 h^{p+1} \|w^{(p+1)}\|_{L_1(\Omega_j)}  \le C_1 h^{p+1} \tau h N^{p+1}_{n,j},
\end{equation}
\begin{equation} \label{EPL2}
\|u^e_j - u^p_j\|_{L_2(\Omega_j)}  \le C_2 h^{p+1} \|w^{(p+1)}\|_{L_2(\Omega_j)} \le C_2 h^{p+1} \tau\sqrt{h}  N^{p+1}_{n,j},
\end{equation}
and
\begin{equation} \label{EPinfinity}
\|u^e_j-u^p_j\|_{L_{\infty}(\Omega_j)}  \le C_3 h^{p+1} \|w^{(p+1)}\|_{L_{\infty}(\Omega_j)} \le C_3 h^{p+1} \tau N^{p+1}_{n,j}.
\end{equation}
Here $C_1, C_2$ and $C_3$ are the projection error constants in the
reference  cell. Consequently, in the whole domain $\Omega$, let
$N^{p+1}_n = \max_j N^{p+1}_{n,j}$, and let
\begin{equation*}
u^p(t,x) = \begin{cases} u^p_L(t,x),  & x \in \Omega^e_L, \\
u^p_R(t,x),  & x \in \Omega^e_R, \\
u^p_j(t,x),  & x \in \Omega_j, \qquad j=0,\cdots,i-2,i+2,\cdots,m-1,
\end{cases}
\end{equation*}
we have
\begin{equation} \label{EPL1G}
\|u^e - u^p\|_{L_1(\Omega)}   \le C_1 |\Omega| h^{p+1} \tau  N^{p+1}_{n},
\end{equation} 
\begin{equation} \label{EPL2G}
\|u^e - u^p\|_{L_2(\Omega)}  \le C_2 \sqrt{|\Omega|}  h^{p+1} \tau  N^{p+1}_{n},
\end{equation} 
and, for the cell boundary terms,
\begin{equation} \label{EPboundary}
\left(\sum_{j=0}^{i-1} + \sum_{j=i+2}^{m-1} |u^e(x_{j+1/2}^-)-u^p(x_{j+1/2}^-)|^2 \right)^{1/2} \le C_3 \sqrt{|\Omega|} h^{p+1/2}  \tau N^{p+1}_{n}.
\end{equation}
\end{lemma}
\noindent {\bf Proof.} This is a direct application of Lemma \ref{lemmaSS} and the Bramble-Hilbert Lemma. \# \\

Now we are ready to state and prove the estimate of $u^p
- u^h$ , then we will conclude this subsection with the  main
theorem to estimate $\tilde{u}(t_{n+1}) - u^h(t_{n+1})$.

\begin{lemma} \label{lemmaPH}
 Let
\begin{equation*}
\hat{u}^p(t,x) = \begin{cases} u^p_L(t,x),  & x \in \Omega^h_L, \\
u^p_R(t,x),  & x \in \Omega^h_R, \\
u^p_j(t,x),  & x \in \Omega_j, \qquad j=0,\cdots,i-2,i+2,\cdots,m-1.
\end{cases}
\end{equation*}
There is a computable constant $Q_n$, depending on $S^p_n$, such
that, for any $t \in (t_n,t_{n+1}]$,
\begin{equation}  \label{PHL2}
\|\hat{u}^p(t,x) - u^h(t,x)\|_{L^2(\Omega)} \le \tau h^{p+1} Q_n.
\end{equation}
Consequently,
\begin{equation}  \label{PHL1}
\|\hat{u}^p(t,x) - u^h(t,x)\|_{L^1(\Omega)} \le \tau h^{p+1}Q_n \sqrt{|\Omega|} .
\end{equation}
Moreover, (\ref{PHL2}) implies 
\begin{equation} \label{PHinfinity}
\|\hat{u}^p(t,x)-u^h(t,x)\|_{L^{\infty}(\Omega)} \le  \tau h^{p+1/2} \tilde{Q}_n
\end{equation}
for some computable constant $\tilde{Q}_n$. 
\end{lemma}
\noindent {\bf Proof.} The proof of Lemma 4.5 in \cite{Rumsey} can be applied to the two subdomains $\Omega^h_L \cup (\cup_{j=0}^{i-2}\Omega_j)$ and $\Omega^h_R \cup (\cup_{j=i+2}^{m-1}\Omega_j)$. Although the statement of Lemma 4.5 in \cite{Rumsey} only covers (\ref{PHL1}),  the proof does cover (\ref{PHL2}) . By using a standard scaling argument, one can prove (\ref{PHinfinity}) from (\ref{PHL2}).   \# \\

$u^p$ and $\hat{u}^p$ are different only in $(\Omega^e_L \cup \Omega^e_R) \, \diagdown (\Omega^e_L \cap \Omega^h_L) \diagdown (\Omega^e_L \cap \Omega^h_L)$, which is the small interval between the semi-discrete shock $x_{sh}(t)$ and the auxiliary shock $x_{se}(t)$. One more lemma is needed to estimate the distance between $x_{sh}(t)$ and  $x_{se}(t)$.
\begin{lemma} \label{shockEH}
There is a computable function $C_s(\cdot)$, such that, for  $t \in (t_n, t_{n+1}]$,
$$
| x_{se}(t) - x_{sh}(t) | \le \tau h^{p+1} C_s(S^p_n).
$$
\end{lemma}
\noindent {\bf Proof.} Since 
$$
\frac{d x_{se}}{dt} = \frac{f(u^e_L(x_{se}))-f(u^e_R(x_{se}))}{u^e_L(x_{se})-u^e_R(x_{se})},
$$
and
$$
\frac{d x_{sh}}{dt} = \frac{f(u^h(x_{sh}^-))-f(u^h(x_{sh}^+))}{u^h(x_{sh}^-)-u^h(x_{sh}^+)},
$$
notice $x_{se}(t_n) = x_{sh}(t_n) = x_{sc,n}$, we get
\begin{eqnarray*}
| x_{se}(t) - x_{sh}(t) | &\le&  \int_{t_n}^t \left| \frac{f(u^e_L(x_{se}))-f(u^e_R(x_{se}))}{u^e_L(x_{se})-u^e_R(x_{se})} 
-\frac{f(u^h(x_{sh}^-))-f(u^h(x_{sh}^+))}{u^h(x_{sh}^-)-u^h(x_{sh}^+)} \right| dt \\
&\le&  \hspace{0.1in} \int_{t_n}^t  \left| \frac{f(u^e_L(x_{se}))-f(u^e_R(x_{se}))}{u^e_L(x_{se})-u^e_R(x_{se})} 
-  \frac{f(u^e_L(x_{sh}))-f(u^e_R(x_{sh}))}{u^e_L(x_{sh})-u^e_R(x_{sh})}  \right| dt \\
&&  + \int_{t_n}^t  \left| \frac{f(u^e_L(x_{sh}))-f(u^e_R(x_{sh}))}{u^e_L(x_{sh})-u^e_R(x_{sh})} 
-\frac{f(u^h(x_{sh}^-))-f(u^h(x_{sh}^+))}{u^h(x_{sh}^-)-u^h(x_{sh}^+)} \right| dt .
\end{eqnarray*}
Here $u^e_L(x_{sh})$ and $u^e_R(x_{sh})$ are well defined, because of the intersection $\Delta ABC$ in the previous figure. Now it is elementary to derive that
\begin{eqnarray*}
&&  \left| \frac{f(u^e_L(x_{se}))-f(u^e_R(x_{se}))}{u^e_L(x_{se})-u^e_R(x_{se})} 
 -  \frac{f(u^e_L(x_{sh}))-f(u^e_R(x_{sh}))}{u^e_L(x_{sh})-u^e_R(x_{sh})}  \right|     \\
&& \quad \le \frac{2\beta}{|u^e_L(x_{se})-u^e_R(x_{se})|} \left[ |u^e_L(x_{se})-u^e_L(x_{sh})| + |u^e_R(x_{se})-u^e_R(x_{sh})| \right] \\
&& \quad \le  \frac{2\beta (N^1_{n,i-1} + N^1_{n,s})}{|u^e_L(x_{se})-u^e_R(x_{se})|} \,\, |x_{se} - x_{sh}| ,
\end{eqnarray*}
and
\begin{eqnarray*}
&& \left| \frac{f(u^e_L(x_{sh}))-f(u^e_R(x_{sh}))}{u^e_L(x_{sh})-u^e_R(x_{sh})} 
-\frac{f(u^h(x_{sh}^-))-f(u^h(x_{sh}^+))}{u^h(x_{sh}^-)-u^h(x_{sh}^+)} \right|    \\
&& \quad \le \frac{2\beta}{|u^e_L(x_{sh})-u^e_R(x_{sh})|} \left[ |u^e_L(x_{sh})-u^h_L(x_{sh}^-)| + |u^e_R(x_{sh})-u^h_R(x_{sh}^+)| \right] \\
&& \quad \le \frac{2\beta (G_L + G_R)}{|u^e_L(x_{sh})-u^e_R(x_{sh})|} \,\, \tau h^{p+1/2}.
\end{eqnarray*}
In the last inequality we have used (\ref{EPinfinity}) of Lemma \ref{lemmaEP} and (\ref{PHinfinity}) of Lemma \ref{lemmaPH}. $G_L$ and $G_R$ can be computed from the coefficients over there. Now, since the shock height is bounded away from zero for a fully developed shock, by combining the last a few inequalities, we can get a computable constant $C_s$, such that
\begin{equation}
| x_{se}(t) - x_{sh}(t) | 
\le \tau^2 h^{p+1/2}  C_s  .
\end{equation}
By tracing the proof, one can see the dependence of $C_s$  on $S^p_n$. \# \\

Combining Lemma \ref{lemmaTE}, Lemma \ref{lemmaEP}, Lemma \ref{lemmaPH} and Lemma \ref{shockEH}, we have the following theorem.
\begin{theorem} \label{spaceTheorem}
There is a computable constant $\F(S^p_n)$, depending on the flux
function $f$, the known constants of the interpolation/projection
error estimates, the known constants of the inverse inequalities,
and the components of the spatial smoothness indicator $S^p_n$, such
that
\begin{equation} \label{errorSemi}
\|\tilde{u}(t_{n+1})  - u^h(t_{n+1})\|_{L_1(\Omega)} \le \tau
h^{p+1} \F(S^p_n).
\end{equation}
\end{theorem}
\noindent {\bf Proof.} By using the sequence of lemmas of this subsection in the subdomains
$$\Omega_{left} = (\Omega^h_L \cap \Omega^e_L) \cup \overline{\cup_{j=0}^{i-2}\Omega_j},
$$ 
$$\Omega_{right} = (\Omega^h_R \cap \Omega^e_R) \cup \overline{\cup_{j=i+2}^{m-1}\Omega_j},
$$
and
$$
\Omega_{shock} = \text{ the interior of } ( \, \Omega \, \diagdown \Omega_{left} \diagdown \Omega_{right} \, )
$$
appropriately, one can easily see that (\ref{errorSemi}) is valid.
\#\\

\subsection{Estimating $u^h(t_{n+1}) - \hat{u}^c_{n+1}$, proof of
(\ref{timeError}) }

The temporal smoothness indicator $T^k_n$ informs us about the
boundedness of the temporal derivatives of $x_{sh}$ and $u^h$
at $t=t_n$. However, what we need to make sure is that the boundedness of $T^k_n$
implies the boundedness of the derivatives of $x_{sh}(t)$ and 
$u^h(t,x)$ for all $t \in [t_n,t_{n+1}]$.

\begin{lemma} \label{lemmaHS}
There is a computable constant K, depending on the spatial
smoothness indicator $S^p_n$, such that, for all $t \in
[t_n,t_{n+1}]$,
\begin{equation} \label{uh}
\| u^h(t) \| _{L_{\infty}(\Omega)} \le \| u^c_n \|
_{L_{\infty}(\Omega)} + K h.
\end{equation}
For each integer $l \in \{1,\cdots, k+1\}$, there is a triple of
computable constants $(b_l, c_l, d_l)$, depending on $S^p_n$, such that, for all $t \in
[t_n,t_{n+1}]$,
\begin{equation} \label{shockSpeedBound}
\left| \frac{d^l x_{sh}}{dt^l} \right| \le b_l
\end{equation}
and
\begin{equation} \label{uhl}
\left\| \frac{\partial^l}{\partial t ^l} u^h(t) \right\| _{L^{\infty}(\Omega)}
\le (1+c_l h^{\alpha}) \left\| \frac{\partial^l}{\partial t ^l} u^h( t_n
)\right\| _{L^{\infty}(\Omega)} + d_l h^{\alpha}.
\end{equation}
To be more specific, for each $l \in \{2,\cdots, k+1\}$, $(b_l, c_l, d_l)$ depends on its predecessors.
\end{lemma}
\noindent {\bf Proof.}
First, one can prove the boundedness of $\|u^h(t,x)\|_{L^{\infty}(\Omega)}$ in $[t_n,t_{n+1}]$ as in the proof of Lemma 4.7 in \cite{Rumsey}. With this done, the boudedness of $\left( \frac{d^l x_{sh}}{dt^l}, \frac{\partial^l}{\partial t ^l} u^h(t) \right)$ can be proven inductively for $l=1,\cdots, k+1$.  

In the case of $l=1$, it is easy to see from (\ref{semiRankine}) that $\left| \frac{d x_{sh}}{dt} \right| \le \beta$. As for the boundedness of $\left\| \frac{\partial}{\partial t } u^h(t,x) \right\| _{L_{\infty}(\Omega)}$, we actually have the proof of Lemma 4.7 in \cite{Rumsey} applying to  the two subdomains left and right of the moving boundary $x_{sh}$. While the proof in \cite{Rumsey} is based on (\ref{secondD}) for regular cells, it can be easily generalized to (\ref{secondL}) for $\Omega^h_L$ and (\ref{secondR}) for $\Omega^h_R$.

After finishing the proof for the case $l=1$, we can see the boundedness of $\frac{d^2 x_{sh}}{dt^2}$ from (\ref{secondS}). Then the boundedness of $\frac{\partial^2}{\partial t^2} u^h(t,x)$ can be proven in a similar way. Inductively, one can do the proof for $l=3,\cdots,k+1.$
 \# \\

Lemma \ref{lemmaHS} has guaranteed the boundedness of the temporal derivatives of the semi-discrete solution $(x_{sh}(t),u^h(t,x))$ in $[t_n,t_{n+1}]$. Hence a standard error analysis on a stiff ODE system results in the next theorem.

\begin{theorem} \label{timeTheorem}
There is a computable function $\G(T^k_n,S_n^{p})$, such that
\begin{equation} \label{errorHC}
\|u^h(t_{n+1}) - \hat{u}^c_{n+1}\|_{L_1(\Omega)}
\le \tau_n^{k+1} \G(T^k_n,S_n^{p}) .
\end{equation}
\end{theorem}
\noindent {\bf Proof.} Because $x_{sc,n+1}$ and $u^c_{n+1}(x)$ are computed by the Kunga-Kutta scheme, the ODE error analysis should give us
$$
| \, x_{sh}(t_{n+1}) - x_{sc,n+1} \, | \le \tau^{k+1} C_1
$$
and 
$$
\| u^h(t_{n+1},x) - \hat{u}^c_{n+1}(x) \|_{L^{\infty}(\Omega_j)} \le\tau^{k+1} C_2, \qquad \Omega_j=\Omega_0,\cdots,\Omega_{i-2},\Omega^h_L,\Omega^h_R,\Omega_{i+2},\cdots,\Omega_{m-1}
$$
for some computable constants $C_1$ and $C_2$ depending on the result of Lemma \ref{lemmaHS}. Now (\ref{errorHC}) follows right away. \#\\

\subsection{The error resulting from the $\Omega_R$ transition}

Every time the two special cells are redome, there is no error resulting from converting the old $\Omega_L$ to $\Omega_{i-1}$ and the new $\Omega_L$. However, the $L^2$-projection from the union $\Omega_R \cup \Omega_{i+2}$ to the new $\Omega_R$ causes some error. This error can be easily estimated. Denote the projection by $P$, and the result of the projection by $P \hat{u}^c_{n+1}$. Then we have the following result.

\begin{lemma} \label{projectionErrorOmegaR} 
There is a computable function $\mathcal{H}(\cdot)$, such that, for the old $\Omega_R = (x_{sc,n+1}, x_{i+3/2})$ and $\Omega_{i+2} = (x_{i+3/2},x_{i+5/2})$,
$$
 \|u^c_{n+1}(x) - \hat{u}^c_{n+1}(x)\|_{L^1(\Omega_R \cup \Omega_{i+2})}  \le h^{p+2} \mathcal{H}(\hat{S}^p_{n+1}),
$$
hence
$$
 \|{u}^c_{n+1}(x) - \hat{u}^c_{n+1}(x)\|_{L^1(\Omega)}  \le h^{p+2} \mathcal{H}(\hat{S}^p_{n+1}) ,
$$
where $\hat{S}^p_{n+1}$ only depend on $\hat{u}^c_{n+1}(x)$ in the union $\Omega_R \cup \Omega_{i+2}$ before the transition.
\end{lemma}
\noindent {\bf Proof.} Let
$$ 
\hat{M}^l_{n+1,i+2} = \frac{\partial^l}{\partial x ^l} \hat{u}^c_{n+1}(x_{i+3/2}^+), \qquad
\hat{L}^l_{n+1,i+2} = \frac{\partial^l}{\partial x ^l} \hat{u}^c_{n+1}(x_{i+3/2}^-),
$$
$$
\hat{A}^l_{n+1,i+2} = \frac{ \hat{M}^l_{n+1,i+2} + \hat{L}^l_{n+1,i+2}}{2}, \qquad \hat{J}^l_{n+1,i+2} = \frac{ \hat{M}^l_{n+1,i+2} - \hat{L}^l_{n+1,i+2}}{2}, 
$$
and
$$
v(x) = \sum_{l=0}^p \frac{\hat{A}^l_{n+1,i+2}}{l!} (x-x_{i+3/2})^l.
$$
Then, since the jumps should satisfy the numerical smoothness condition $|\hat{J}^l_{n+1,i+2}| = \hat{D}^l_{n+1,i+2} h^{p+2-l(1+\alpha)}$ with $\hat{D}^l_{n+1,i+2} = \O(1)$, and it can be verified in the computation, we have
\begin{eqnarray*}
&& \|v(x)-\hat{u}^c_{n+1}(x)\|^2_{L^2(\Omega_R \cup \Omega_{i+2})} \\
&&  = \left\| \sum_{l=0}^p \frac{\hat{J}^l_{n+1,i+2}}{2 \, l!} (x-x_{i+3/2})^l \right\|^2_{L^2(\Omega_{i+2}) }
+   \left\| \sum_{l=0}^p \frac{-\hat{J}^l_{n+1,i+2}}{2 \, l!} (x-x_{i+3/2})^l \right\|^2_{L^2(\Omega_R) } 
 \le h^{2p+3} \hat{C}^2,
\end{eqnarray*}
for a computable $\hat{C}$ depending on the smoothness indicator $\hat{D}^l_{n+1,i+2}$ and $\hat{M}^l_{n+1,i+2}$. Now,
\begin{eqnarray*}
&&  \|u^c_{n+1}(x) - \hat{u}^c_{n+1}(x)\|_{L^1(\Omega_R \cup \Omega_{i+2})} \\
&& = \|P\hat{u}^c_{n+1}(x) - \hat{u}^c_{n+1}(x)\|_{L^1(\Omega_R \cup \Omega_{i+2})} \\
&&  \le \sqrt{2h} \,\, \|P\hat{u}^c_{n+1}(x) - \hat{u}^c_{n+1}(x)\|_{L^2(\Omega_R \cup \Omega_{i+2})} \\
&& \le \sqrt{2h} \,\, \|v(x) - \hat{u}^c_{n+1}(x)\|_{L^2(\Omega_R \cup \Omega_{i+2})} \\
&& \le h^{p+2} \sqrt{2} \hat{C}.
\end{eqnarray*}
Taking $\H(\hat{S}^p_{n+1}) = \sqrt{2} \hat{C}$,  we have proven the lemma. \# \\

{\bf Remark.} This projection error is added to the global error, not for every time step, but only every time the shock location passes a regular cell boundary point. Therefore, in the power $h^{p+2}$, one of the power is to be lost to the error accumulation. The total contribution of these projections to the global error is $\mathcal{O}(h^{p+1})$. \\

Combining Theorem \ref{spaceTheorem}, Theorem \ref{timeTheorem} and Lemma \ref{projectionErrorOmegaR}, we have proven  Theorem \ref{main}.

\section{Numerical experiments}

\begin{figure}
  \includegraphics[width=6.2in,height=3.6in]{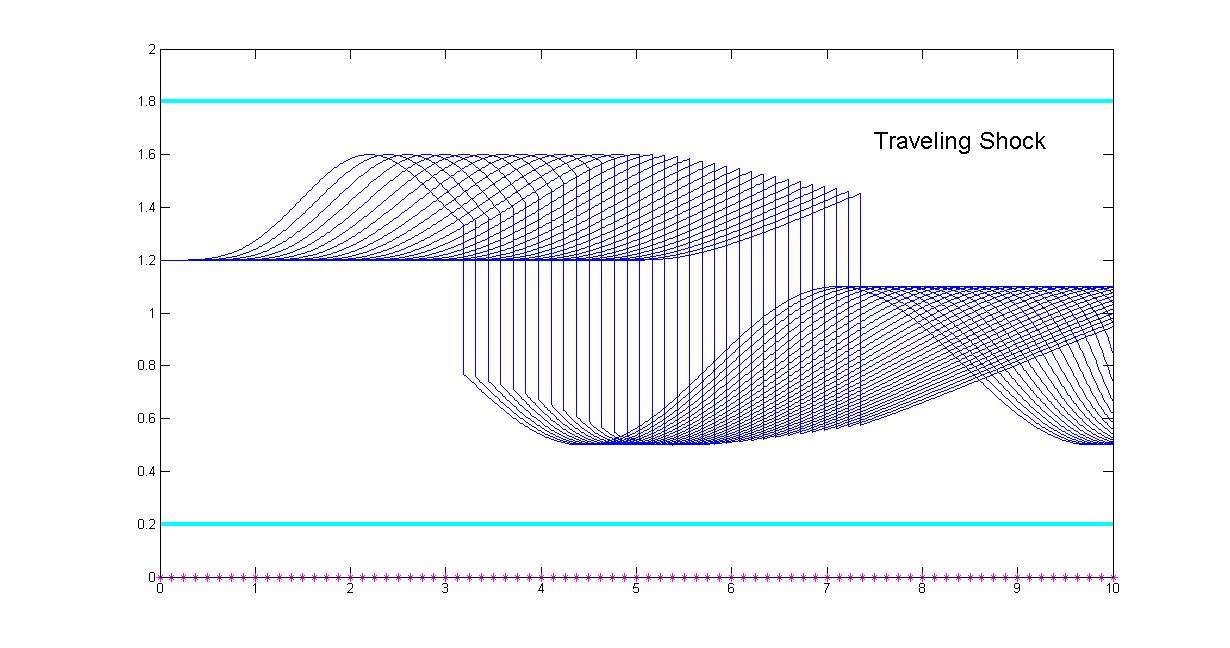}
  \caption{The traveling shock, computed for $t \in [0,4]$}\label{travelingShock}
\end{figure}

As mentioned in Section 2, we solve the Burgers' equation $u_t = \left( \frac{u^2}{2} \right)_x$ as the numerical example. The initial condition with a shock is plotted in Figure \ref{initialShock}, where $(a,b) = (0,10)$, $x_s(0) = c = 3.18$, the boundary condition is $u_a(t)=1.2$, and the initial condition is
\begin{equation*}
u_I(x) = \begin{cases} 1.2 + 0.4 \, \sin^4 (x/1.4),              &    \text{if $ x \in (0,x_s(0))$}, \\
                                      0.8 - 0.3 \, \sin((x-3.1)/0.85),          &  \text{if $ x \in (x_s(0), 10)$}. \end{cases}
\end{equation*}
In Figure \ref{travelingShock}, we plot the solution at different times, to show that the shock and the rest of the wave are computed with neither smearing nor overshooting. Moreover, one can see the solution's further sharpening and rarefaction, the change of the wave speed depending on $u$, the growth of the shock height, the extrema
of the smooth pieces kept constant until swallowed by the shock, the subsequent total variation diminishing, etc.. All these suggest that the computed solution is of high quality.

In order to show some more informative data about the convergence rate of the solution, we will show the computed  smoothness indicators and the convergence rate of the shock in the next two subsections. In the third subsection, we will show a case where, due to an inappropriate choice of the time step size, the scheme makes the numerical solution non-smooth numerically. 

\subsection{The numerical smoothness indicators}

\begin{figure}
  \includegraphics[width=6.0in,height=3.5in]{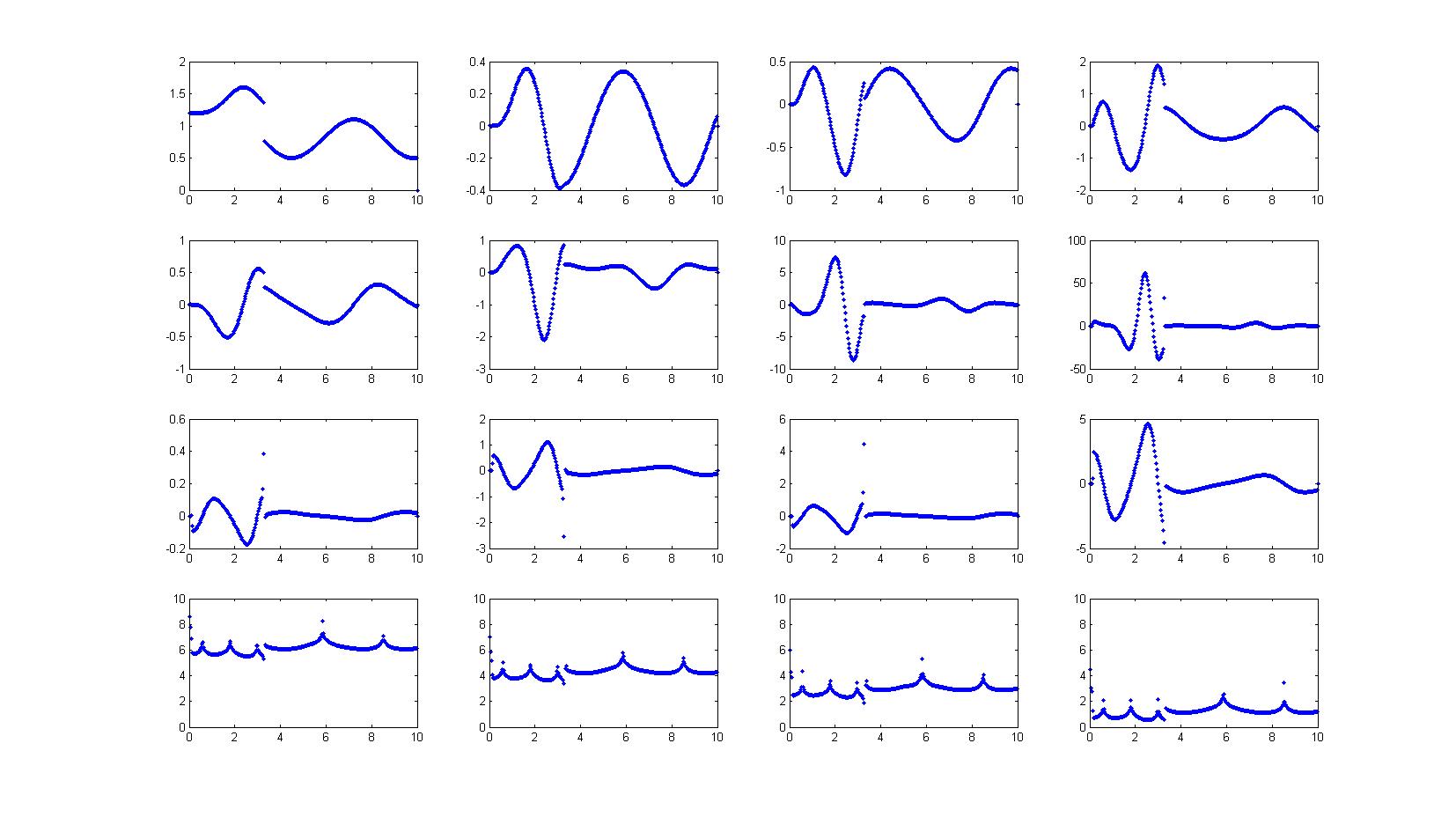}
  \caption{The components of the smoothness indicator, n=71}\label{indicator71}
\end{figure}

Table \ref{indicatorContents} shows the contents of the pictures of the numerical smoothness indicators. Figure \ref{indicator71} shows, as an example, the values of the  smoothness indicators after time step $n=71$ for the case $h=1/32$ and $\tau = \frac{h}{20}$. While the first three rows of a smoothness indicator figure show the contents of the smoothness indicator as defined before, the last row shows the orders of the jumps $J^l_{n,j}$ of the derivatives, because 
$$
\log_h J^l_{n,j} = \log_h [ D^0_{n,j} h^{p+2-l(1+\alpha)} ] =  p+2-l(1+\alpha) + \log_h D^0_{n,j} .
$$
If we indeed have $D^0_{n,j} = \O(1)$, then $\log_h J^0_{n,j} \approx  p+2-l(1+\alpha)$, where $\alpha=1/p$. It is interesting to observe the shape similarity of the four logarithm curves, which deliver some information about the solution more explicitly.

For the purpose of supporting the lemmas, theorems and proofs of this paper, it suffices to show the boundedness of the smoothness indicator. These bounded values enter the error estimates in the lemmas and theorems. However, the set of indicator values deliver extremely abundant information about the solution, including the quality of the solution, the discontinuities of different orders, the influence of the cell sizes and its changes, the need for a local refinement if it occurs, the anti-smoothness behavior of the scheme if it occurs, and etc.. A separate technical report will be written to expand the numerical results.

\begin{table}
\caption{The contents of the numerical smoothness indicator figures }
\centering
 \label{indicatorContents}
\begin{tabular}{|p{1in}|p{1in}|p{1in}|p{1in}|}
\hline
$M^0_{n,j}$		&    $M^1_{n,j}$	&   $M^2_{n,j}$	&   $M^3_{n,j}$		\\ \hline
$u^h_t(t_n,x)$ 	&    $u^h_{tt}(t_n,x)$	&   $u^h_{ttt}(t_n,x)$	&   $u^h_{tttt}(t_n,x)$		\\ \hline
$D^0_{n,j}$ 	&    $D^1_{n,j}$	&   $D^2_{n,j}$	&   $D^3_{n,j}$		\\ \hline
$\log_h J^0_{n,j}$ &$\log_hJ^1_{n,j}$&$\log_hJ^2_{n,j}$&$\log_h J^3_{n,j}$ 	\\ \hline
\end{tabular}
\end{table}

\subsection{The convergence rate of the shock location}

According to the main theorem \ref{main}, the shock location should converge at the rate $\mathcal{O}(h^4+\tau^3)$ when $p=3$ and $k=3$. Since the error resulting from space discretization and time discretization cannot be separated, we should make $\tau = \mathcal{O}(h^{1+1/3})$ in order to balance the error. This is also consistent with the assumption $\tau \le \gamma h^{1+\alpha}$ in the theorems and proofs. Since we do not have the true solution available for this example,  instead of the differences between a numerical solution and the true solution, we use the differences between the consecutive numerical solutions in the converging sequence. 

In Table \ref{T1}, we have the computed shock locations for six different cell sizes $h=\frac{1}{2},\frac{1}{4}, \frac{1}{8}, \frac{1}{16}, \frac{1}{32}$ and $\frac{1}{64}$,  where the time step sizes are decided so that $\tau \approx \frac{h^{1+1/3}}{10}$. The last column of the table clearly  indicates that the shock location is indeed converging at the rate of $h^4$.

\begin{table} 
\caption{Convergence rate of the shock,  when $\tau \approx \frac{h^{1+1/3}}{10}$ }
\centering
 \label{T1}
\begin{tabular}{|p{0.5in}|p{0.75in}|p{1.33in}|p{1.28in}|p{1.15in}|}
\hline
$h$ & $\tau \approx \frac{h^{1+1/3}}{10}$ & $x_{sc,n}(h), \, t=n\tau=4.0$ & $|x_{sc,n}(h)-x_{sc,n}(h/2)|$ & $\frac{|x_{sc,n}(h)-x_{sc,n}(h/2)|}{|x_{sc,n}(h/2)-x_{sc,n}(h/4)|}$ \\ \hline
1/2 	&  1/25		&  7.353087041875537	&  3.80697 e-05	&  29.001 	\\ \hline
1/4 	&  1/64 	&  7.353125111551957	&  1.31272 e-06	&  18.744 	\\ \hline
1/8 	&  1/160	&  7.353126424273150 	&  7.00341 e-08	&  16.986	\\ \hline
1/16	&  1/400 	&  7.353126494307232	&  4.12303 e-09	&  17.529	\\ \hline
1/32	&  1/1024 	&  7.353126498430261	&  2.35212 e-10	&		\\ \hline
1/64	&  1/2560 	&  7.353126498665473	&  			&		\\ \hline
\end{tabular}
\end{table}

In Table \ref{T2}, we show the results obtained by using the same set of cell sizes, but a slightly different $\gamma$. Namely,  $\tau \approx \frac{h^{1+1/3}}{8}$. Again, we observe a convergence rate of $h^4$.

\begin{table} 
\caption{Convergence rate of the shock, when $\tau \approx \frac{h^{1+1/3}}{8}$}
\centering
\label{T2}
\noindent \begin{tabular}{|p{0.5in}|p{0.75in}|p{1.33in}|p{1.28in}|p{1.15in}|}
\hline
$h$ & $\tau \approx \frac{h^{1+1/3}}{8}$ & $x_{sc,n}(h), \, t=n\tau=4.0$ & $|x_{sc,n}(h)-x_{sc,n}(h/2)|$ & $\frac{|x_{sc,n}(h)-x_{sc,n}(h/2)|}{|x_{sc,n}(h/2)-x_{sc,n}(h/4)|}$ \\ \hline
1/2 	&  1/20 	&  7.353059324500445 	&  6.40705 e-05 	&  21.716 	\\ \hline
1/4 	&  1/50		&  7.353123395074651 	&  2.95044 e-06	&  20.421 	\\ \hline
1/8 	&  1/128 	&  7.353126345518764 	&  1.44479 e-07       	&  17.224	\\ \hline
1/16	&  1/320 	&  7.353126489997938 	&  8.15233 e-09        	&  16.290 	\\ \hline
1/32	&  1/800 	&  7.353126498150260	&  5.00450 e-10	&		\\ \hline
1/64	&  1/2048 	&  7.353126498650710	&  			&		\\ \hline
\end{tabular}
\end{table}

In Table \ref{T3}, the time step size $\tau$ is determined by  $\tau \approx \frac{h^{1+1/6}}{10}$. Here, we seem to observe a convergence rate around $h^{3.5}$, certainly less than optimal in space.

\begin{table} 
\caption{Convergence rate of the shock,  when $\tau \approx \frac{h^{1+1/6}}{10}$}
\centering
\label{T3}
\noindent \begin{tabular}{|p{0.5in}|p{0.75in}|p{1.33in}|p{1.28in}|p{1.15in}|}
\hline
$h$ & $\tau \approx \frac{h^{1+1/6}}{10}$ & $x_{sc,n}(h), \, t=n\tau=4.0$ & $|x_{sc,n}(h)-x_{sc,n}(h/2)|$ & $\frac{|x_{sc,n}(h)-x_{sc,n}(h/2)|}{|x_{sc,n}(h/2)-x_{sc,n}(h/4)|}$ \\ \hline
1/2 	&  1/20 	&  7.353059324500445 	&  6.36534 e-05	&  19.392 	\\ \hline
1/4 	&  1/48 	&  7.353122977913548	&  3.28251 e-06        	&  14.866 	\\ \hline
1/8 	&  1/112 	&  7.353126260423104	&  2.20802 e-07       	&  13.792	\\ \hline
1/16	&  1/256 	&  7.353126481224717	&  1.60092 e-08	&  12.131	\\ \hline
1/32	&  1/576	&  7.353126497233903	&  1.31970 e-09	&		\\ \hline
1/64	&  1/1280 	&  7.353126498553604	&  			&		\\ \hline
\end{tabular}
\end{table}

In Table \ref{T4}, we show the results obtained by using  $\tau = \frac{h}{16}$. With such time step sizes, we observe a convergence rate of $h^3$, losing a whole order.

\begin{table} 
\caption{Convergence rate of the shock,  when $\tau = \frac{h}{16}$}
\centering
\label{T4}
\noindent \begin{tabular}{|p{0.5in}|p{0.75in}|p{1.33in}|p{1.28in}|p{1.15in}|}
\hline
$h$ & $\tau=h/16$ & $x_{sc,n}(h), \, t=n\tau=4.0$ & $|x_{sc,n}(h)-x_{sc,n}(h/2)|$ & $\frac{|x_{sc,n}(h)-x_{sc,n}(h/2)|}{|x_{sc,n}(h/2)-x_{sc,n}(h/4)|}$ \\ \hline
1/2 	&  1/32 	&  7.353107562791714 	&  1.74588 e-05 	&  14.149 	\\ \hline
1/4 	&  1/64 	&  7.353125111551957 	&  1.23397 e-06	&  9.0930 	\\ \hline
1/8 	&  1/128 	&  7.353126345518764 	&  1.35706 e-07 	&  8.8269	\\ \hline
1/16	&  1/256 	&  7.353126481224717	&  1.53742 e-09	&  8.4015 	\\ \hline
1/32	&  1/512 	&  7.353126496598886	&  1.82994 e-09	&		\\ \hline
1/64	&  1/1024 	&  7.353126498428824	&  			&		\\ \hline
\end{tabular}
\end{table}

\subsection{An anti-smoothing phenomenon}

It is well-known that, for the RKDG method, $\tau$ needs to be taken fairly small, otherwise the scheme can be ``unstable". However, the strengthened CFL condition $\tau \le \gamma h^{1+\alpha}$ is not necessary according to \cite{qz2009}. On the other side, if a linear CFL condition $\tau \le \tilde{\gamma} h$ is assumed, it is not known what $\tilde{\gamma}$ should be. The following numerical experiment shows that, for this numerical example, $\tilde{\gamma} = \frac{1}{10}$ is not sufficiently small to guarantee ``numerical stability". 

\begin{figure}
  \includegraphics[width=6.0in,height=3.5in]{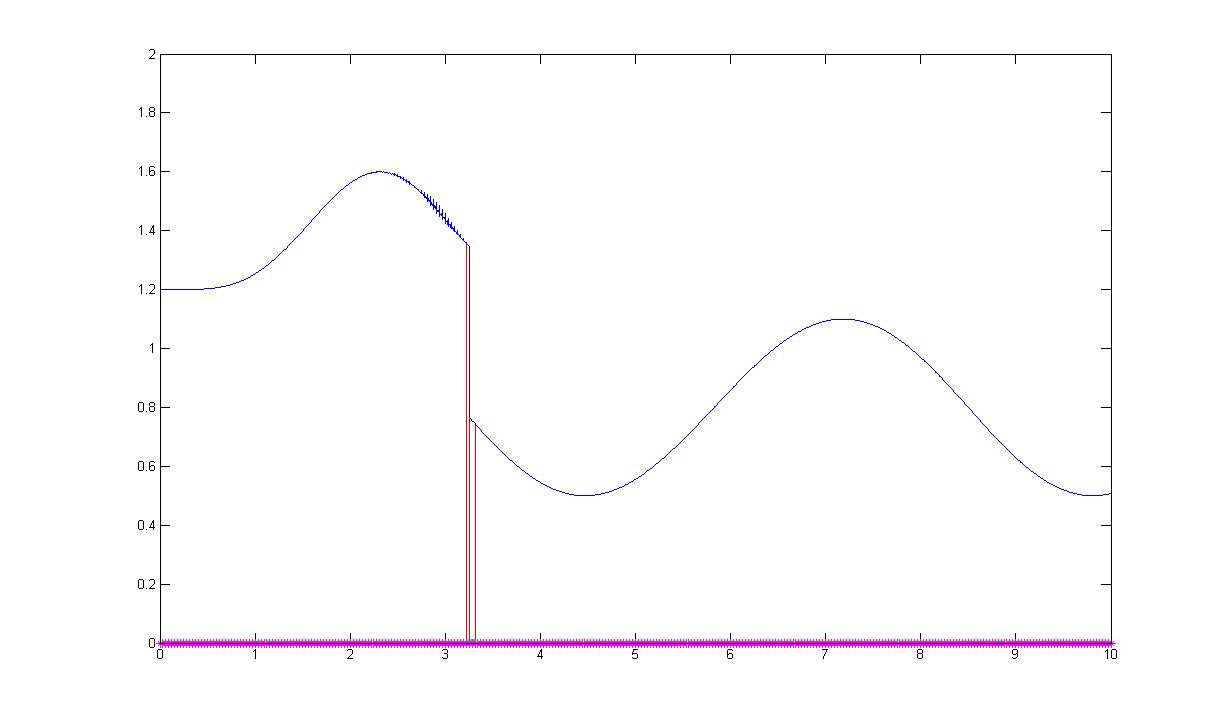}
  \caption{Loss of numerical smoothness, n=24}\label{antiSmooth}
\end{figure}

\begin{figure}
  \includegraphics[width=6.0in,height=3.5in]{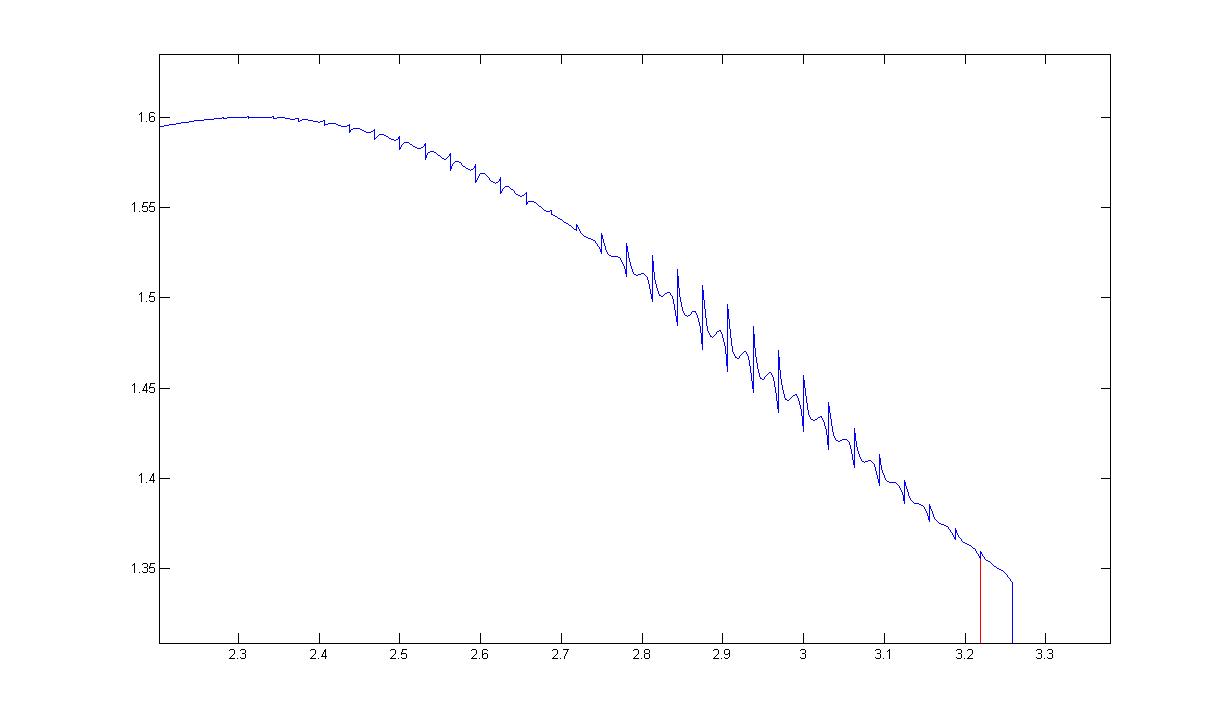}
  \caption{Zoom in at the loss of smoothness, n=24}\label{zoom24}
\end{figure}

\begin{figure}
  \includegraphics[width=6.0in,height=3.5in]{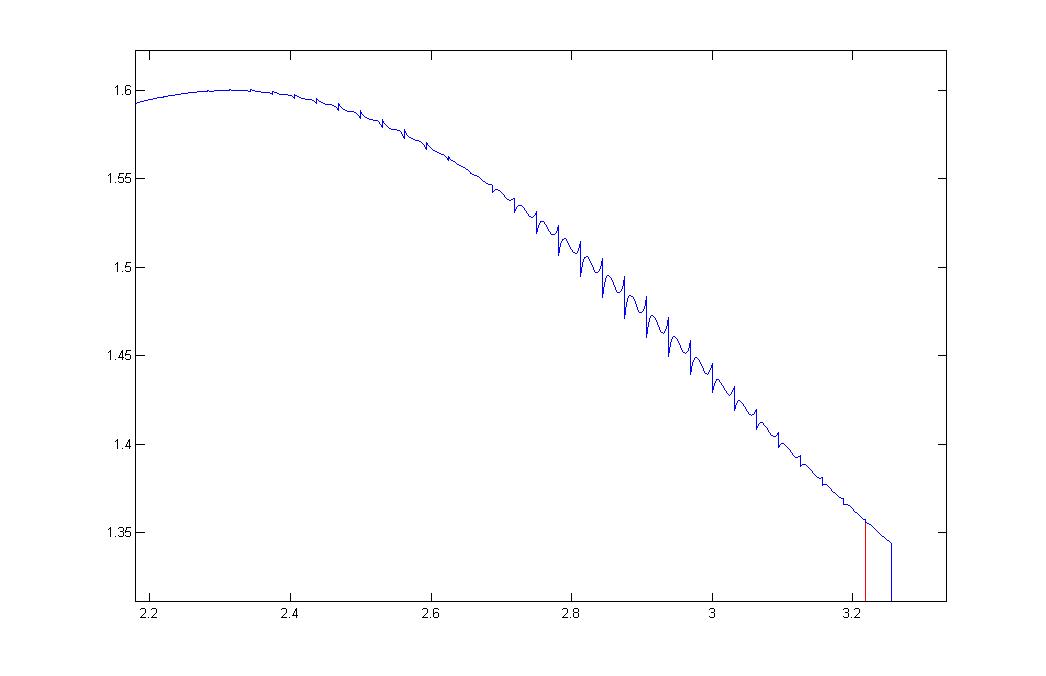}
  \caption{Zoom in at the loss of smoothness, n=23}\label{zoom23}
\end{figure}

We tried $\tau = \frac{h}{10}$ for $h=\frac{1}{2}, \frac{1}{4}, \frac{1}{8}, \frac{1}{16},$ and $\frac{1}{32}$. While the first four cases seem OK, ``numerical instability" occurred in the case of $h=\frac{1}{32}$, $\tau = \frac{1}{320}$. Figure \ref{antiSmooth} shows that the numerical solution became oscillatory in the supposedly smooth piece of the solution left of the shock. This is at the time $t=\frac{24}{320}$, where $n=24$. Since it is in the upwind side, the oscillation has nothing to do with the shock, because the Godunov flux in the scheme does not deliver information in the upwind direction at all. To have a better look, by a zooming-in to the oscillatory piece of the solution, Figure \ref{zoom24} shows the details of the oscillation, where the cubic polynomials in each cell have become really cubic-looking. To trace the source of the oscillation, Figure \ref{zoom23} shows the details of the oscillation one time steps earlier.

% Figure \ref{zoom20} and Figure \ref{zoom19} show the details of the oscillation a few time steps earlier. 

A few things about the oscillation can be mentioned. (1) The size of the oscillation is growing over the steps. (2) The third derivative in each cell is changing sign from step to step. (3) The size of the oscillation is changing gradually from cell to cell.   Our guess is that the TVDRK-3 scheme has just stepped out of its ``stability region"  in the process of $h = 10 \tau \rightarrow 0$, although it is too hard to prove anything due to the nonlinearity. Here, TVD-stability is certainly lost, but Lax-stability seems still alright. Nevertheless, the optimal convergence rate has been lost, because numerical smoothness is obviously lost. In fact, the numerical smoothness indicator has detected the loss of smoothness even earlier.  

\begin{figure}
  \includegraphics[width=6.0in,height=3.5in]{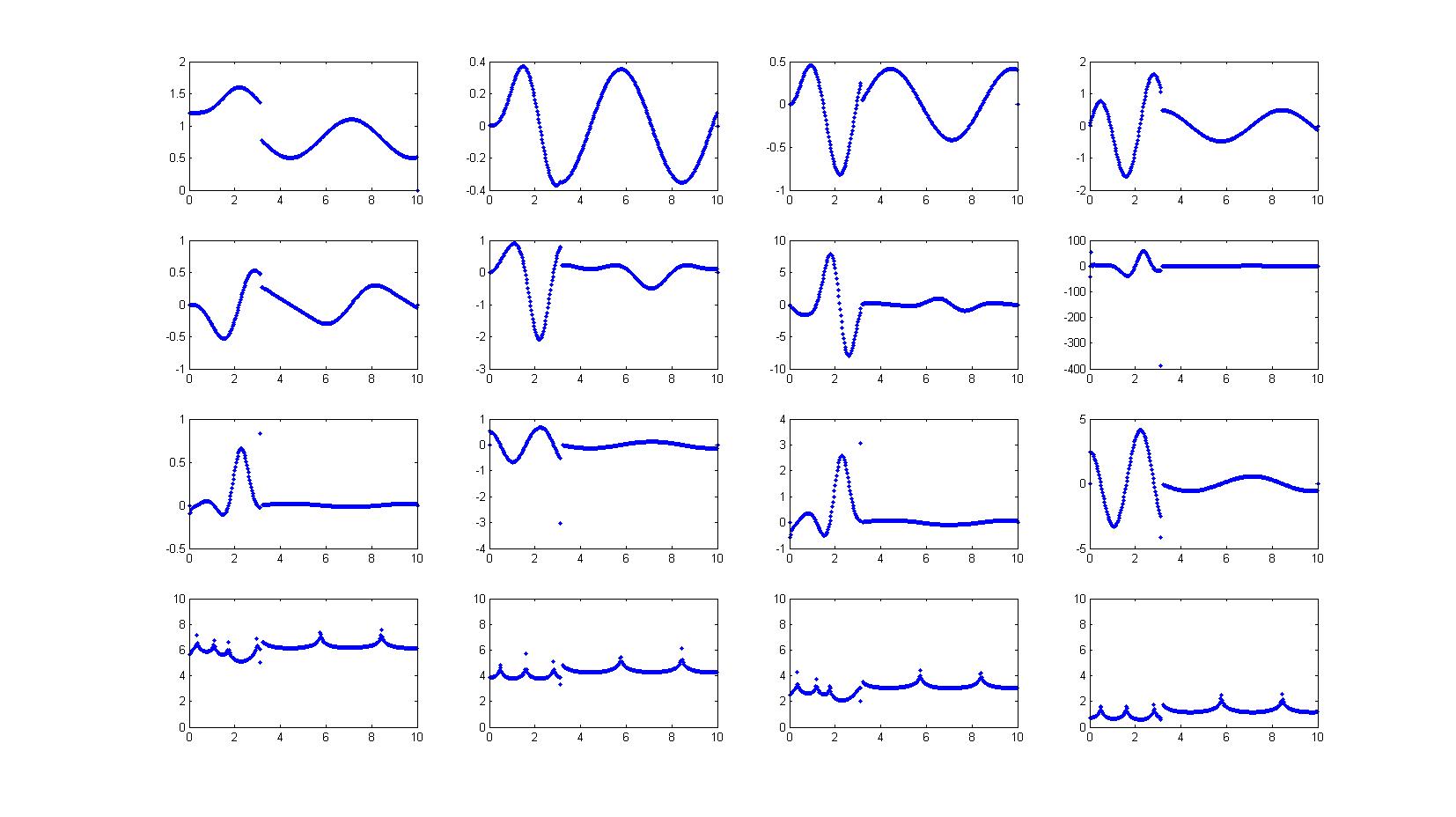}
  \caption{$D^l_{n,j}$ has grown in size, n=2}\label{instability2}
\end{figure}

\begin{figure}
  \includegraphics[width=6.0in,height=3.5in]{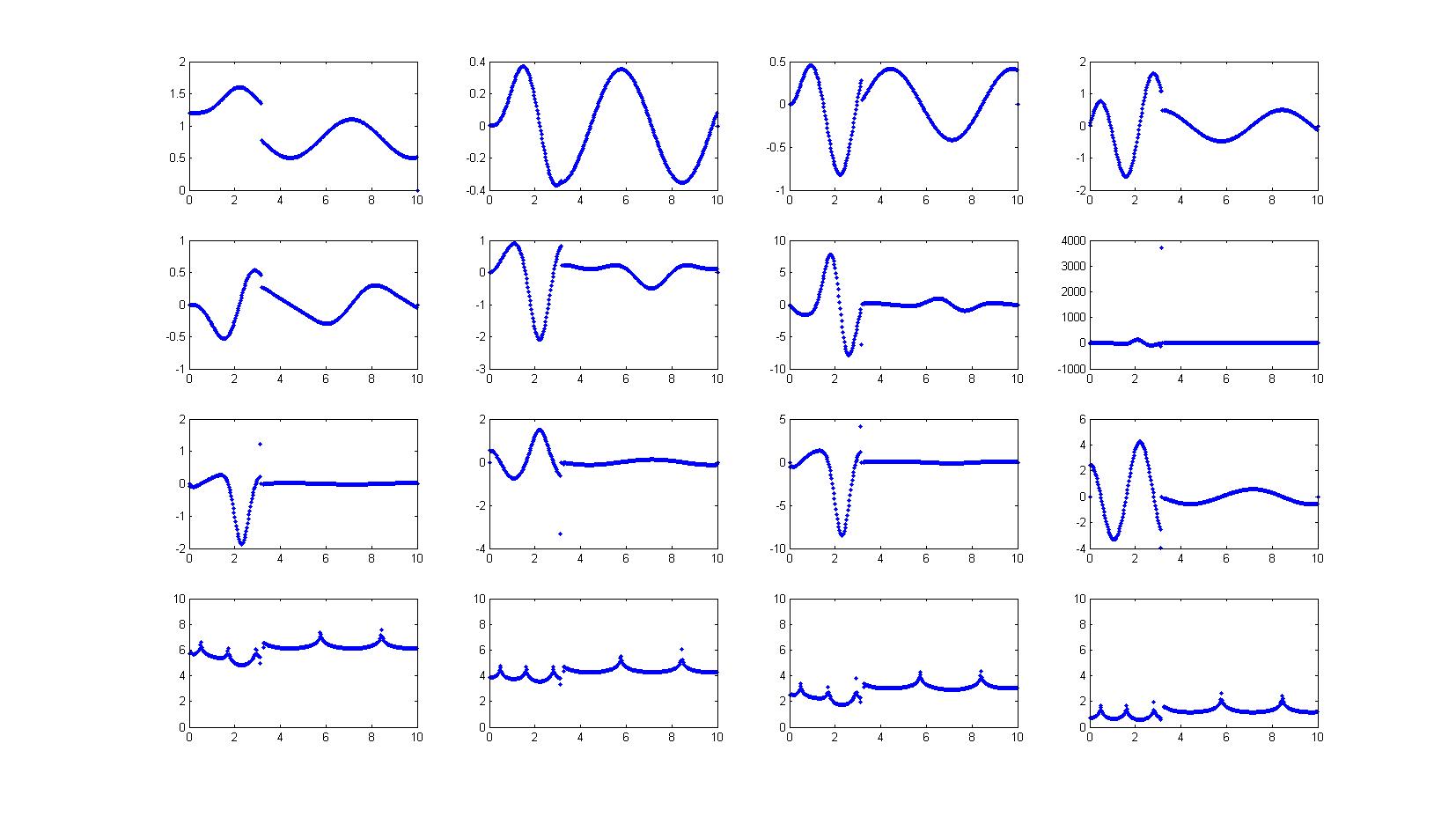}
  \caption{$D^l_{n,j}$ has grown in size and changed sign, n=2}\label{instability3}
\end{figure}

Figure \ref{instability2} and Figure \ref{instability3} show the numerical smoothness indicators after time step 2 and 3. Focusing in the region before the shock, one can see that, not only the sizes of  $D^l_{n,j}$ were growing fast, their signs changed between the two steps. This is a typical early symptom of the loss of numerical smoothness, also understood as loss of numerical stability \cite{strikwerda}. It is worth to mention that, for the loss to be visible from the solution itself, it will take 15 to 20 more time steps. The numerical smoothness indicator works as a microscope, detecting the trouble when it is still  ``invisible". At the beginning the trouble is only showing through $D^0_{n,j}$, but all the other high order terms of the smoothness indicator also show problems a few steps later. In fact, at $n=3$ there has been almost no damage down to the solution itself. By switching to $\tau = h/20$ staring from $n=4$, we still obtained a high quality global solution. 

%\section{Conclusion remarks}

%\noindent {\bf A.} 

\bibliographystyle{amsplain}

%\end{article}

\end{document}